\documentclass[12pt]{article}
\usepackage[latin1]{inputenc}
\usepackage{amssymb}
\usepackage{hyperref}
\usepackage{amsmath}
\usepackage{latexsym}
\usepackage{times}
\usepackage{cite}
\usepackage{amssymb}
\usepackage{amsfonts}
\usepackage{amscd}
\usepackage{xcolor}
\usepackage{amstext,amsmath,amssymb,amsfonts}
\newtheorem{theorem}{Theorem}

\newtheorem{corollary}[theorem]{Corollary}

\newtheorem{definition}[theorem]{Definition}

\newtheorem{lemma}[theorem]{Lemma}

\newtheorem{proposition}[theorem]{Proposition}
\newtheorem{remark}[theorem]{Remark}

\newcommand{\beq}{\begin{eqnarray}}
\newcommand{\eeq}{\end{eqnarray}}
\newcommand{\beqs}{\begin{eqnarray*}}
	\newcommand{\eeqs}{\end{eqnarray*}}
\newcommand{\bpro}{\begin{pro}}
	\newcommand{\epro}{\end{pro}}
\newcommand{\blem}{\begin{lem}}
	\newcommand{\elem}{\end{lem}}
\newcommand{\bdfn}{\begin{dfn}}
	\newcommand{\edfn}{\end{dfn}}
\newcommand{\bcor}{\begin{cor}}
	\newcommand{\ecor}{\end{cor}}
\newcommand{\bthm}{\begin{thm}}
	\newcommand{\ethm}{\end{thm}}
\newcommand{\bex}{\begin{ex}}
	\newcommand{\eex}{\end{ex}}
\newcommand{\brmk}{\begin{rmk}}
	\newcommand{\ermk}{\end{rmk}}
\newcommand{\bpr}{\begin{pr}}
	\newcommand{\epr}{\end{pr}}
\newcommand{\benum}{\begin{enumerate}}
	\newcommand{\eenum}{\end{enumerate}}
\newcommand{\bitem}{\begin{itemize}}
	\newcommand{\eitem}{\end{itemize}}

\newcommand{\cqfd}{\hfill{\square}}
\chardef\bslash=`\\
\numberwithin{equation}{section}
\numberwithin{table}{section}
\numberwithin{theorem}{section}

\begin{document}
	\begin{flushright}
		ICMPA-MPA/2019
	\end{flushright}
	\begin{center}
		{\Large { $\mathcal{R}(p,q)-$ deformed  combinatorics: full characterization and illustration}}\\
		\vspace{0,5cm}
		Mahouton Norbert Hounkonnou$^{*}$ and Fridolin Melong\\
		\vspace{0.5cm}
		{\em International Chair in Mathematical Physics
			and Applications}
		{\em (ICMPA-UNESCO Chair), }
		{\em University of Abomey-Calavi,}
		{\em 072 B.P. 50 Cotonou, Republic of Benin}\\	
	\end{center}
	\today
	
	\vspace{0.5 cm}
	{\it Abstract}
	
	This paper addresses a theory of  $\mathcal{R}(p,q)$-deformed combinatorics in discrete probability. It mainly focuses on  $\mathcal{R}(p,q)$-deformed   factorials, binomial coefficients, Vandermonde's formula, Cauchy's formula, binomial and negative binomial formulae, factorial and binomial moments, and Stirling numbers. Moreover, the $\mathcal{R}(p,q)-$ Stirling numbers of the second kind and the $\mathcal{R}(p,q)-$ Bell numbers  for graphs are also derived. Related  relevant properties are investigated and discussed. Finally,  as a concrete illustration,  the developed   formalism is displayed for the well known  generalized $q-$ Quesne deformed quantum algebra to construct the corresponding 
	deformed combinatorics, as a particular case.
	
		\tableofcontents
	\section{Introduction}
	Combinatorial theory is a major branch of  mathematics, which has  applications in many fields such as computer science (languages, graphs, intelligent computing), natural and social sciences, biomedicine, molecular biology, operational research, engineering, and business\cite{L,M,R}.
	Combinatorial theory  and discrete mathematical methods play an important role, and occupy a central position in the theory of discrete probability. The most prominent of these methods are the combinatorial enumerative methods and the basic methods of finite difference computation. A considerable number of stochastic experiments or phenomena   in discrete probability theory can be described by the stochastic models of distributions \cite{CA4}.
	
	Chung and Kang developed a new combinatorics called $q-$ combinatorics and investigated the significance of $q-$ permutations and $q-$ combinations\cite{CK}. The idea of
	$q-$ analogs can be traced back to Euler in the 1700's who studied $q-$ series, especially specializations of theta functions. 
	Meanwhile, 
	in \cite{CA1}, Ch. A. Charalambides examined basic $q-$ combinatorics  and $q-$ hypergeometric series. The $q-$power,  $q-$factorial,  
	$q-$binomial coefficient of a real number and two $q-$ Vandermonde's ($q-$factorial convolution)
	formulae were derived. The $q-$ analogs of the Cauchy's formulas were also  investigated in  \cite{A}. 
	
	Furthermore, the $q-$Stirling
	numbers of the first and second kinds, which are the coefficients of the expansions
	of $q-$factorials into $q-$ powers, and of $q-$powers into $q-$factorials, respectively,   were presented. Moreover, the Stirling numbers of the second kind and their generalizations were studied by several authors, (see for instance \cite{K1, K2, K3} and references theiren). The $q-$ Stirling numbers of the second kind and $q-$ Bell numbers for graphs were also  analyzed in \cite{BS}. Corcino and Barientos \cite{RB} established many properties for  $q-$ analogs of  Stirling numbers. The vertical and horizontal recursion relations and the generating function were also computed. Besides,  two parameters $(p,q)-$ Stirling numbers, which are  generating functions for the joint distribution of pair statistics, were described in \cite{WW}. In addition, a theory of  $(p,q)-$ analogs of  binomial coefficients was elaborated. Some
	properties and identities like  the triangular, vertical, and  horizontal recursion relations,
	the generating function,  the orthogonality and inverse relations were derived in \cite{RC}. 
	
	Later in 2010,  $\mathcal{R}(p,q)-$ deformed quantum  algebra was introduced by Hounkonnou and Bukweli \cite{HB1} as a generalization of known deformed quantum algebras. The same authors also performed $\mathcal{R}(p,q)-$ differentiation and   integration, and deducted all relevant  particular cases of $q-$ and $(p,q)-$ deformations \cite{HB}. This opens a novel route for developing  the theory of $\mathcal{R}(p,q)-$ analogs of special numbers and combinatorics.
	
	This paper provides a general formalism, which enables an easier construction  of  a 
	combinatorial theory from  deformed quantum algebras existing in the literature by assigning concrete suitable expressions to the function $\mathcal{R}$ and related specific meromorphic functions of the theory.	
	Especially, $\mathcal{R}(p,q)-$ analogs of  factorials,  binomial coefficients,  Vandermonde's formula, Cauchy's formula,  binomial formula and  negative binomial formula, factorial moments,  binomial moments, Stirling numbers and  Bell numbers on a graph are investigated and discussed. Furthermore, the  case of the generalized $q-$ Quesne deformed combinatorics is derived to illustrate the presented formalism. From this new generalization, developing  deformed combinatorial theories associated with  other particular cases of deformed quantum algebras,   known and spread in the literature, appears as a matter of triviality.
	
	The  paper is organized as follows. Section $2$ is devoted to basic notations, definitions and results  related to $\mathcal{R}(p,q)-$ quantum algebras and  calculus, and to basic $q-$ combinatorics. In Section $3,$  the fundamentals  of $\mathcal{R}(p,q)-$ deformed combinatorics and the derivation of   relevant properties are exposed. 
	The generalized $q-$ Quesne combinatorics  is derived as a case study. Some concluding remarks are addressed in Section $4.$ 
	
	\section{Preliminaries}
	In this section, we briefly recall the main definitions, notations and known results used in the sequel. For more details, the reader can refer to \cite{CA1,CA2, CA3, CA4, RC}, \cite{HB, HB1} and \cite{J1,J1a}. 
	
	The coherent states introduced by Quesne\cite{Q} can be associated with the
	$q-$deformed algebra satisfying the relations
	\begin{gather}\label{p20}
	[N, a^{\dagger}]= a^{\dagger}\mbox{,}\quad [N, a]= -a,\nonumber\\
	a\,a^{\dagger} - a^{\dagger}\,a =q^{-N-1}\quad\mbox{or}\quad q\,a\,a^{\dagger} - a^{\dagger}\,a =\mathbf{1}
	\end{gather}
	where $0<q<1,$ and the Quesne number is defined by:
	\begin{equation*}\label{p22}
	[u]^{Q}= {1-q^{-u}\over q-1}.
	\end{equation*}
	The Quesne algebra is a particular case of the Kalnins-Miller-Mukherjee algebra \cite{KMM}
	with $\ell= 1, \lambda= 0.$ Furthermore, Hounkonnou and Ngompe Nkouankam \cite{Hounkonnou&Ngompe07a}    generalized the $q-$ Quesne algebra with the generators satisfying the relations
	\begin{gather}\label{p23}
	[N, a^{\dagger}]= a^{\dagger}\mbox{,}\quad [N, a]= -a,\nonumber\\
	p^{-1}\,a\,a^{\dagger} - a^{\dagger}\,a =q^{-N-1}\quad\mbox{or}\quad q\,a\,a^{\dagger} - a^{\dagger}\,a =p^{N+1}
	\end{gather}
	where $0<q <p \leq 1.$ Their generalized $q-$ Quesne number is given as follows:
	\begin{equation*}\label{p25}
	[u]_{p,q}^{Q}= {p^u-q^{-u}\over q-p^{-1}}.
	\end{equation*}

	Let now $\mathcal{R}$ be a meromorphic function defined on $\mathbb{C}\times\mathbb{C}$ by \begin{equation}\label{r10}
	\mathcal{R}(u,v)= \sum_{s,t=-l}^{\infty}r_{st}u^sv^t,
	\end{equation}
	converging in the complex disc $\mathbb{D}_{R}=\left\lbrace z\in\mathbb{C}/ |z|<R\right\rbrace,$
	where $r_{st}$ are complex numbers, $l\in\mathbb{N}\backslash\left\lbrace 0\right\rbrace,$ and $R$ is the radius of convergence of the series (\ref{r10}).
	Let us consider the set of holomorphic functions $\mathcal{O}(\mathbb{D}_{R})$ defined on $\mathbb{D}_{R}$.
	
	In the sequel, when no possible confusion arises, $p$ and $q$ will designate two real numbers satisfying $0<q<p\leq 1.$
	
	\begin{definition}\cite{HB}
		Let $P$ and $Q$ be two linear operators on $\mathcal{O}(\mathbb{D}_{R})$. Then, for $\varPsi\in\mathcal{O}(\mathbb{D}_{R})$, we have \begin{equation*}\label{r2}
		Q:\varPsi\longmapsto Q\varPsi(z):= \varPsi(qz),
		\end{equation*}
		\begin{equation*}\label{r3}
		P:\varPsi\longmapsto P\varPsi(z):= \varPsi(pz).
		\end{equation*}
		The $(p,q)-$derivative and the $(p,q)-$number are defined, respectively, by \cite{CJ}: \begin{equation*}\label{r4}
		D_{p,q} :\varPsi\longmapsto D_{p,q}\varPsi(z):=\frac{\varPsi(pz)-\varPsi(qz)}{z(p-q)},
		\end{equation*}
		\begin{equation*}
		[n]_{p,q}:=\frac{p^n-q^n}{p-q},
		\end{equation*}
		while the $\mathcal{R}(p,q)-$derivative is given by \cite{HB}: 
		\begin{equation*}\label{r5}
		D_{\mathcal{R}( p,q)}:=D_{p,q}\frac{p-q}{P-Q}\mathcal{R}( P,Q)=\frac{p-q}{p^{P}-q^{Q}}\mathcal{R}(p^{P},q^{Q})D_{p,q}.
		\end{equation*}
	\end{definition}
	\begin{definition}\cite{HB1}
		The  $\mathcal{R}(p,q)-$number and the  $\mathcal{R}(p,q)-$ factorials are defined, respectively, as follows:
		\begin{equation*}
		[n]_{\mathcal{R}(p,q)}:=\mathcal{R}(p^n,q^n)\quad\mbox{for}\quad n\geq 0,
		\end{equation*}
		\begin{equation*}\label{s0}
		\mathcal{R}!(p^\kappa,q^\kappa):=\left \{
		\begin{array}{l}
		1\quad\mbox{for}\quad \kappa=0\\
		\\
		\mathcal{R}(p,q)\cdots\mathcal{R}(p^\kappa,q^\kappa)\quad\mbox{for}\quad \kappa\geq 1,
		\end{array}
		\right .
		\end{equation*}
		and the  $\mathcal{R}(p,q)-$ binomial coefficient is given by:
		\begin{equation*}\label{bc}
		\left[\begin{array}{c} m  \\ \kappa\end{array} \right]_{\mathcal{R}(p,q)} := \frac{\mathcal{R}!(p^m ,q^m)}{\mathcal{R}!(p^\kappa,q^\kappa)\mathcal{R}!(p^{m-\kappa},q^{m-\kappa})},\quad m,\kappa=0,1,2,\cdots,\quad m\geq \kappa.
		\end{equation*}
	\end{definition}
	More details on $\mathcal{R}(p,q)-$ deformed quantum algebras, $\mathcal{R}(p,q)-$differentiation and integration can be found in \cite{HB,HB1}.
	To be complete, let us briefly recall some notions about known $q-$ combinatorics pertaining to our development in the sequel.
	\begin{definition}
		The $q-$ shifted factorial is given by
		\begin{equation}\label{qbc1}
		\big(1\oplus y\big)^n_q:= \prod_{i=0}^{n-1}\big(1 + yq^{i-1}\big) = \sum_{\kappa=0}^{n}\,q^{\kappa \choose 2}\bigg[\begin{array}{c}
		n \atop \\ \kappa
		\end{array}\bigg]_{q}\,y^{\kappa}.
		\end{equation}
	\end{definition}
	The $(p,q)-$ binomial coefficients are given by \cite{RC}:
	\begin{equation}\label{l1}
	\prod_{i=1}^{n}(p^{i-1} + u\,q^{i-1})= \sum_{\kappa=0}^{n}\,p^{{n-\kappa \choose 2}}\,q^{\kappa \choose 2}\bigg[\begin{array}{c}
	n \atop \\ \kappa
	\end{array}\bigg]_{p,q}\,u^{\kappa}.
	\end{equation}
	For  real numbers $x$ and $y,$  the $q-$ analog of the  Cauchy's formula is given
	by \cite{CA1}
	\begin{equation}\label{qc1}
	\bigg[\begin{array}{c} x+y  \\ n\end{array} \bigg]_{q} = \sum_{\kappa=0}^{n-1}q^{\kappa(y-n+\kappa)}\,\bigg[\begin{array}{c} x  \\ \kappa\end{array} \bigg]_{q}\,\,\bigg[\begin{array}{c} y  \\ n-\kappa \end{array} \bigg]_{q},
	\end{equation}	
	or, equivalently,
	\begin{equation}\label{qc2}
	\bigg[\begin{array}{c} x+y  \\ n\end{array} \bigg]_{q} = \sum_{\kappa=0}^{n-1}q^{(n-\kappa)(x-\kappa)}\,\bigg[\begin{array}{c} x  \\ \kappa\end{array} \bigg]_{q}\,\,\bigg[\begin{array}{c} y  \\ n-\kappa \end{array} \bigg]_{q}.
	\end{equation}
	Furthermore,
	the following orthogonality relations hold \cite{CA2}:
	\begin{equation}\label{qbn6}
	\sum_{x=\kappa}^{n}(-1)^{n-x}\,q^{{n-x \choose 2}}\bigg[\begin{array}{c} n   \\ x \end{array} \bigg]_{q}\,\bigg[\begin{array}{c} x   \\ \kappa \end{array} \bigg]_{q}=\delta_{n,\kappa},
	\end{equation}
	and
	\begin{equation}\label{qbn7}
	\sum_{x=\kappa}^{n}(-1)^{x-\kappa}q^{{x-\kappa \choose 2}}\bigg[\begin{array}{c} n   \\ x \end{array} \bigg]_{q}\,\bigg[\begin{array}{c} x   \\ \kappa \end{array} \bigg]_{q}=\delta_{n,\kappa}
	\end{equation}
	where  $\delta_{n,\kappa}$ is the Kronecker delta, $n$ and $\kappa$ are positive integers.	
	The $n^{th}-$ order of a $q-$ factorial number ${[x]_{n,q}}$ is written as a polynomial of the $q-$ number as follows \cite{GW}:
	\begin{equation}\label{qsnf}
	{[x]_{n,q}} = q^{-{n\choose 2}}\sum_{\kappa=0}^{n}s_q(n,\kappa)[x]_q^\kappa,
	\end{equation} 
	or,
	\begin{equation}\label{qsns}
	[x]^n_q = \sum_{\kappa=0}^{n}q^{{\kappa \choose 2}}S_q(n,\kappa){[x]_{\kappa,q}}
	\end{equation}
	where the coefficients $s_q(n,\kappa)$ and $S_q(n,\kappa)$ are called  $q-$ Stirling numbers of the first and second kinds, respectively. 
	
	Let $\kappa$ and $j$ be positive integers. Then, the following relations also hold:
	\begin{equation}\label{qs33}
	{\kappa \choose j} =\sum_{m=j}^{\kappa}(-1)^{m-j}(1-q)^{m-j} s_{q}(m,j)\,\bigg[\begin{array}{c} \kappa  \\ m \end{array} \bigg]_{q}
	\end{equation} 
	and
	\begin{equation}\label{qs33b}
	\bigg[\begin{array}{c} \kappa  \\ j \end{array} \bigg]_{q}=\sum_{m=j}^{\kappa}(-1)^{m-j}(1-q)^{m-j} S_{q}(m,j)\,{\kappa \choose m},
	\end{equation}
	where $x$ is an integer and $0<q<1.$
	Finally, recall that the $q-$deformed probability distribution $g(x)$ of a discrete random variable $X$ is given by Ch. A. Charalambides  \cite{CA1} as
	\begin{equation}\label{fba9}
	g(x) = \sum_{m=x}^{\infty}(-1)^{m-x}q^{{m-x\choose 2}}\bigg[\begin{array}{c} m \\ x \end{array} \bigg]_{q}{\bf E}\Big(\bigg[\begin{array}{c} X \\ m \end{array} \bigg]_{q}\Big),\quad x\in\mathbb{N}
	\end{equation}
	where the series is absolutely convergent, and  ${\bf E}\bigg(\bigg[\begin{array}{c} X \\ m \end{array} \bigg]_{q}\bigg)$ stands for the expectation value of $\bigg[\begin{array}{c} X \\ m \end{array} \bigg]_{q}$.
		\section{$\mathcal{R}(p,q)-$ deformed combinatorics}
		Our aim is to present fundamentals of  a generalization of  the combinatorial theory from the $\mathcal{R}(p,q)-$ deformed quantum  algebra  introduced in \cite{HB1} as a generalization of known deformed quantum algebras.
		We consider
		$\Phi_i(p,q)\in \mathcal{O}(\mathbb{D}_{R}),$ with $i=1, 2,$ 
		depending on the parameters $p$ and $q.$
		\subsection{$\mathcal{R}(p,q)-$ deformed factorials and  binomial coefficients}
		
		\begin{definition}\label{def3b}
			Let $u$ be a real number. Then,	the $\mathcal{R}(p,q)-$ deformed factorial of $u$ of order $\kappa$ is defined by : 
			\begin{equation}\label{02}
			{[u]_{\kappa,\mathcal{R}(p,q)}}
			= \prod_{v=1}^{\kappa}\mathcal{R}(p^{u-v+1},q^{u-v+1}),
			\end{equation}
			where $0<q<p\leq 1$  and $\kappa\in\mathbb{N}\backslash \{0\}. $
		\end{definition}
		The relation (\ref{02}) will also be called the $\kappa^{th}-$ order factorial of the $\mathcal{R}(p,q)-$ deformed number. From the above definition, we derive the following basic property for the $\mathcal{R}(p,q)-$ deformed factorial:
		\begin{equation}\label{fp}
		{[u]_{\kappa+s,\mathcal{R}(p,q)}
			= [u]_{s,\mathcal{R}(p,q)}\,[u-s]_{\kappa,\mathcal{R}(p,q)}}
		,\quad \kappa\in\mathbb{N}\backslash \{0\}\quad\mbox{and}\quad s\in\mathbb{N}\backslash \{0\}.
		\end{equation}
		\begin{lemma}
			Let $x,$ $q,$ and $p$ be real numbers such that $0<q<p\leq 1.$ Then, the $\mathcal{R}(p,q)-$ deformed factorial of $x$ of negative order ${-\kappa}$ is given as follows:
			\begin{equation}\label{08}
			{[x]_{-\kappa,\mathcal{R}(p,q)}:={1\over [x+\kappa]_{\kappa,\mathcal{R}(p,q)}}}
			\mbox{,}\quad \kappa\in\mathbb{N}\backslash\{0\}.
			\end{equation}
			\begin{equation*}\label{06}
			{[x]_{0,\mathcal{R}(p,q)}:=1.}
			\end{equation*}
		\end{lemma}
		{\it Proof:} It uses the relation (\ref{fp}). $\cqfd$
		
		For $x=0,$ the equation (\ref{08}) yields
		\begin{equation*}\label{09}
		{[0]_{-\kappa,\mathcal{R}(p,q)}
			= {1 \over [\kappa]_{\mathcal{R}(p,q)}!}}
		\mbox{,}\quad \kappa \in\mathbb{N}\backslash\{0\}.
		\end{equation*}
		
		\begin{definition}
			The $\mathcal{R}(p,q)-$ deformed binomial coefficient is defined by:
			\begin{equation}\label{010}
			\bigg[\begin{array}{c} u  \\ \kappa \end{array} \bigg]_{\mathcal{R}(p,q)} := 
			{[u]_{\kappa,\mathcal{R}(p,q)}\over[\kappa]_{\mathcal{R}(p,q)}!} 
			\mbox{,}\quad \kappa\in\mathbb{N}\backslash \{0\},\;\; u\in\mathbb{R}.
			\end{equation}
		\end{definition}
		We assume there exist 	$\Phi_i\in \mathcal{O}(\mathbb{D}_{R}),$ with $i=1, 2,$ 
		depending on the parameters $p$ and $q,$ which link 
		the  deformed numbers $\mathcal{R}(p^u,q^u)$ and $\mathcal{R}(p^{-u},q^{-u})-$ as follows:
		\begin{equation}\label{011}
		\mathcal{R}(p^{-u},q^{-u}) := \Big(\Phi_1(p,q)\,\Phi_2(p,q)\Big)^{1-u}\,\mathcal{R}(p^u,q^u).
		\end{equation}
		Then,   the following relations hold: 
		\begin{equation}\label{012}
		\mathcal{R}(p^{u+v},q^{u+v}) = \Phi^{v}_1(p,q)\,\mathcal{R}(p^{u},q^{u}) + \Phi^{u}_2(p,q)\,\mathcal{R}(p^{v},q^{v})
		\end{equation}
		and  
		\begin{equation}\label{013}
		\mathcal{R}(p^{u-v},q^{u-v}) = \Phi^{-v}_1(p,q)\,\mathcal{R}(p^{u},q^{u}) - \Phi^{-v}_1(p,q)\, \Phi^{u-v}_2(p,q)\,\mathcal{R}(p^{v},q^{v}),
		\end{equation}
		$u$ and $v$ being  real numbers.
		\begin{proposition}
			Let $u$ be a natural number and $\kappa$ a positive integer. Then, 
			\begin{equation}\label{014}
			{[u]_{\kappa,\mathcal{R}(p^{-1},q^{-1})}}
			= \Big(\Phi_1(p,q)\,\Phi_2(p,q)\Big)^{-u\kappa + {\kappa +1 \choose 2}}\,{[u]_{\kappa,\mathcal{R}(p,q)}}
			\end{equation}
			\begin{equation}\label{015}
			\mathcal{R}!(p^{-u},q^{-u})= \Big(\Phi_1(p,q)\,\Phi_2(p,q)\Big)^{- {u  \choose 2}}\,\,\mathcal{R}!(p^{u},q^{u})
			\end{equation}
			and
			\begin{equation}\label{016}
			\bigg[\begin{array}{c} u  \\ \kappa \end{array} \bigg]_{\mathcal{R}(p^{-1},q^{-1})} =\Big(\Phi_1(p,q)\,\Phi_2(p,q)\Big)^{-\kappa(u-\kappa)}\,\,\bigg[\begin{array}{c} u  \\ \kappa \end{array} \bigg]_{\mathcal{R}(p,q)},\quad\kappa\in\mathbb{N}.
			\end{equation}
		\end{proposition}
		{\it Proof:} 
		Using the relation (\ref{011}), we get
		\begin{eqnarray}
		{[u]_{\kappa,\mathcal{R}(p^{-1},q^{-1})}}
		&=& \prod_{j=1}^{\kappa}\mathcal{R}(p^{-u+j-1},q^{-u+j-1})\nonumber\\
		&=& \Big(\Phi_1(p,q)\,\Phi_2(p,q)\Big)^{-u\kappa + {\kappa +1 \choose 2}}\,{[u]_{\kappa,\mathcal{R}(p,q)}}.
		\end{eqnarray}
		Furthermore,\begin{eqnarray}
		\mathcal{R}!(p^{-u},q^{-u}) &=& \prod_{j=1}^{u}\mathcal{R}(p^{-j},q^{-j})\nonumber\\
		&=& \prod_{j=1}^{u}\Big(\Phi_1(p,q)\,\Phi_2(p,q)\Big)^{1-j}\mathcal{R}(p^j,q^j)\nonumber\\
		&=& \Big(\Phi_1(p,q)\,\Phi_2(p,q)\Big)^{-{\kappa  \choose 2}}\mathcal{R}!(p^{u},q^{u}).
		\end{eqnarray} 
		Finally, using the relation (\ref{010}), we obtain (\ref{016}).
		$\cqfd$
		\begin{remark}
			The formula (\ref{016}) may be expressed as
			\begin{equation}\label{a016}
			\bigg[\begin{array}{c} u  \\ \kappa \end{array} \bigg]_{\mathcal{R}({q\over p})} =\Phi^{-\kappa(u-\kappa)}_1(p,q)\,\,\bigg[\begin{array}{c} u  \\ \kappa \end{array} \bigg]_{\mathcal{R}(p,q)},\quad\kappa\in\mathbb{N}.
			\end{equation}
			For $\mathcal{R}(u,v)=(p-q)^{-1}(u-v),$ we recover the \textbf{Jagannathan-Srinivassa} binomial coefficient  as:
			\begin{equation*}\label{Ja016}
			\bigg[\begin{array}{c} u  \\ \kappa \end{array} \bigg]_{{q\over p}} =p^{-\kappa(u-\kappa)}\,\,\bigg[\begin{array}{c} u  \\ \kappa \end{array} \bigg]_{p,q},\quad\kappa\in\mathbb{N}.
			\end{equation*}
			The \textbf{generalized $q-$Quesne} binomial coefficient can be obtained by putting $\mathcal{R}(x,y)=(q-p^{-1})^{-1}(x-y^{-1})$ as:
			\begin{equation*}\label{ga016}
			\bigg[\begin{array}{c} u  \\ \kappa \end{array} \bigg]^Q_{{1\over p\,q}} =p^{-\kappa(u-\kappa)}\,\,\bigg[\begin{array}{c} u  \\ \kappa \end{array} \bigg]^Q_{p,q},\quad\kappa\in\mathbb{N}.
			\end{equation*}
		\end{remark}
		\begin{theorem}
			Let $x,$ $p$ and $q$ be  real numbers such that $ 0 < q < p \leq 1,$ and $\kappa$ be a positive integer. Then, the $\mathcal{R}(p,q)-$ deformed binomial coefficients satisfy the following recursion relation:  
			\begin{equation}\label{bc1}
			\bigg[\begin{array}{c} x  \\ \kappa \end{array} \bigg]_{\mathcal{R}(p,q)} = \Phi^{\kappa}_1(p,q)\,\bigg[\begin{array}{c} x -1 \\ \kappa \end{array} \bigg]_{\mathcal{R}(p,q)} + \Phi^{x-\kappa}_2(p,q)\,\bigg[\begin{array}{c} x-1  \\ \kappa-1 \end{array} \bigg]_{\mathcal{R}(p,q)},\quad\kappa\in\mathbb{N}
			\end{equation}
			or,  equivalently,
			\begin{equation}\label{bc2}
			\bigg[\begin{array}{c} x  \\ \kappa \end{array} \bigg]_{\mathcal{R}(p,q)} = \Phi^{\kappa}_2(p,q)\,\bigg[\begin{array}{c} x -1 \\ \kappa \end{array} \bigg]_{\mathcal{R}(p,q)} + \Phi^{x-\kappa}_1(p,q)\,\bigg[\begin{array}{c} x-1  \\ \kappa-1 \end{array} \bigg]_{\mathcal{R}(p,q)},\quad\kappa\in\mathbb{N}
			\end{equation}
			with the initial conditon $\bigg[\begin{array}{c} x  \\ 0 \end{array} \bigg]_{\mathcal{R}(p,q)} :=1.$
		\end{theorem}
		{\it Proof:} {Since
			$[x]_{\kappa,\mathcal{R}(p,q)}=[x]_{\mathcal{R}(p,q)}\,[x-1]_{\kappa-1,\mathcal{R}(p,q)},$    
			$[x-1]_{\kappa,\mathcal{R}(p,q)}=[x-1]_{\kappa-1,\mathcal{R}(p,q)}\,[x-\kappa]_{\mathcal{R}(p,q)},$ using the relations (\ref{012}),  the $\mathcal{R}(p,q)-$ factorial of x}
		satisfies the recursion relation:
		\begin{equation*}
		{[x]_{\kappa,\mathcal{R}(p,q)}=\Phi^{\kappa}_1(p,q)\,[x-1]_{\kappa,\mathcal{R}(p,q)}+\Phi^{x-\kappa}_2(p,q)\,[\kappa]_{\mathcal{R}(p,q)}\,[x-1]_{\kappa-1,\mathcal{R}(p,q)}}
		\end{equation*}
		with condition $[x]^0_{\mathcal{R}(p,q)}:=1.$
		From the expression (\ref{010}), 
		the relation (\ref{bc1}) is deduced. Furthermore,  the expression (\ref{02})  satisfies
		\begin{equation*}
		{[x]_{\kappa,\mathcal{R}(p,q)}=\Phi^{x-\kappa}_1(p,q)\,[\kappa]_{\mathcal{R}(p,q)}\,[x-1]_{\kappa-1,\mathcal{R}(p,q)} + \Phi^{\kappa}_2(p,q)\,[x-1]_{\kappa,\mathcal{R}(p,q)}.}
		\end{equation*}
		Dividing the members of the above equation by $[\kappa]_{\mathcal{R}(p,q)}!$, we obtain (\ref{bc2}), and the proof is achieved. $\cqfd$
		\begin{corollary}
			Let $x,$ and $r$  be real numbers. Then, the following relation holds:
			\begin{equation*}
			\prod_{j=1}^{n}\big(\Phi^{x-r-i+1}_1(p,q) - \Phi^{x-r-i+1}_2(p,q)\big) = \big(\Phi_1(p,q) - \Phi_2(p,q)\big)^n\,{[x-r]_{n,\mathcal{R}(p,q)}}.
			\end{equation*}
		\end{corollary}
		{\it Proof:} It is straightforward by computation. $\cqfd$
		
		Taking $\mathcal{R}(u,v)= v,$ we recuperate a simpler relation under the form
		\begin{equation*}
		\prod_{j=1}^{n}(1 - q^{x-r-j+1}) = (1 - q)^n\,{[x-r]_{n,q}}
		\end{equation*}
		where $0<q< 1.$
		\begin{remark}
			\begin{enumerate}
				\item[(1)] Note that the results obtained by Hounkonnou and Bukweli in \cite{HB} can be retrieved by taking $\mathcal{R}(s,t)=((p^{-1}-q)t)(st-1).
				$
				
				\item [(2)]	The \textbf{generalized $q-$Quesne} formulae are given as follows, with real numbers $p$ and $q$  such that $0<q<p\leq 1$ and $\kappa\in\mathbb{N}\backslash\{0\}:$
				\begin{itemize}
					\item[(i)] 
					
					\begin{equation*}\label{02g}
					{[x]^Q_{\kappa,p,q}}
					= \prod_{v=1}^{\kappa}[x-v+1]^Q_{p,q}.
					\end{equation*}
					\item[(ii)]
					\begin{equation*}\label{010g}
					\bigg[\begin{array}{c} x  \\ \kappa \end{array} \bigg]^Q_{p,q} = {
						{[x]^Q_{\kappa,p,q}}\over [\kappa]^Q_{p,q}!}.
					\end{equation*}
					\item[(iii)]
					\begin{equation*}\label{011g}
					[x]^Q_{p^{-1},q^{-1}} = \Big({p\over q}\Big)^{-1-x}\,[x]^Q_{p,q}.
					\end{equation*}
					\item[(iv)]
					\begin{equation*}\label{014g}
					{[x]^Q_{\kappa,p^{-1},q^{-1}}}
					= \Bigg({p\over q}\Bigg)^{-2\kappa -x\kappa + {\kappa +1 \choose 2}}\,\,{[x]^Q_{\kappa,p,q}}
					,
					\end{equation*}
					\begin{equation*}\label{015g}
					[x]^Q_{p^{-1},q^{-1}}! = \Bigg({p\over q}\Bigg)^{-x - {x+1  \choose 2}}\,\,[x]^Q_{p,q}!,
					\end{equation*}
					and
					\begin{equation*}\label{016g}
					\Bigg[\begin{array}{c} x  \\ \kappa \end{array} \Bigg]^Q_{p^{-1},q^{-1}} =\Bigg({p\over q}\Bigg)^{-\kappa(x-\kappa)}\,\,\Bigg[\begin{array}{c} x  \\ \kappa \end{array} \Bigg]^Q_{p,q}.
					\end{equation*}
				\end{itemize}
			\end{enumerate}
		\end{remark}
		\subsection{$\mathcal{R}(p,q)-$ deformed Vandermonde's and   Cauchy's formulae}
		The  $\mathcal{R}(p,q)-$ deformed Vandermonde's formula, also called $\mathcal{R}(p,q)-$ deformed factorial convolution, is contained in the theorem below, where, for the expression simplification,  we set $\Phi_1(p,q)=\epsilon_1$ and $\Phi_2(p,q)=\epsilon_2;$ $x,$ $y,$ $p,$ and $q$ are  real numbers with $ 0 < q < p \leq 1.$
		\begin{theorem}
			The $\mathcal{R}(p,q)-$ deformed Vandermonde's formula is given by:
			\begin{equation}\label{vd1}
			{[x+y]_{n,\mathcal{R}(p,q)}}
			= \sum_{\kappa=0}^{n}\bigg[\begin{array}{c} n  \\ \kappa\end{array} \bigg]_{\mathcal{R}(p,q)}\epsilon_1^{\kappa(y-n+\kappa)}\epsilon_2^{(n-\kappa)(x-\kappa)}{[x]_{\kappa,\mathcal{R}(p,q)}[y]_{n-\kappa,\mathcal{R}(p,q)}}
			\end{equation}
			or, equivalently,
			\begin{equation}\label{vd2}
			{[x+y]_{n,\mathcal{R}(p,q)}}
			= \sum_{\kappa=0}^{n}\bigg[\begin{array}{c} n  \\ \kappa\end{array} \bigg]_{\mathcal{R}(p,q)}\epsilon_1^{(n-\kappa)(x-\kappa)}\epsilon_2^{\kappa(y-n+\kappa)}{[x]_{\kappa,\mathcal{R}(p,q)}[y]_{n-\kappa,\mathcal{R}(p,q)}},
			\end{equation}
			where $n$ is a positive integer.
		\end{theorem}
		{\it Proof:} For $n\in \mathbb{N}\backslash \{0\},$ we consider the following expression:
		\begin{equation}
		T_n(x;y)_{\mathcal{R}(p,q)}:= \sum_{\kappa=0}^{n}\bigg[\begin{array}{c} n  \\ \kappa\end{array} \bigg]_{\mathcal{R}(p,q)}\epsilon_1^{\kappa(y-n+\kappa)}\epsilon_2^{(n-\kappa)(x-\kappa)}{[x]_{\kappa,\mathcal{R}(p,q)}[y]_{n-\kappa,\mathcal{R}(p,q)}}.
		\end{equation}
		For $n=1,$ we have
		$T_1(x;y)_{\mathcal{R}(p,q)}
		= \mathcal{R}(p^{x+y},q^{x+y}).$ Using the recursion relation (\ref{bc1}) and 
		\begin{eqnarray}\label{vdc}
		\big[x+y-n+1\big]_{\mathcal{R}(p,q)}
		{[x]_{\kappa,\mathcal{R}(p,q)}[y]_{n-\kappa-1,\mathcal{R}(p,q)}}
		&=& \epsilon_1^{y-n+\kappa+1}
		{[x]_{\kappa+1,\mathcal{R}(p,q)}[y]_{n-\kappa-1,\mathcal{R}(p,q)}}
		\nonumber\\ &+& \epsilon_2^{x-\kappa}{[x]_{\kappa,\mathcal{R}(p,q)}[y]_{n-\kappa,\mathcal{R}(p,q)}},	
		\end{eqnarray}
		we obtain
		\begin{eqnarray*}\label{spa}
			T_n(x;y)_{\mathcal{R}(p,q)} &=& \sum_{\kappa=0}^{n-1}\bigg[\begin{array}{c} n-1  \\ \kappa\end{array} \bigg]_{\mathcal{R}(p,q)}\epsilon_1^{\kappa(y-n+\kappa+1)}\epsilon_2^{(n-\kappa)(x-\kappa)}{[x]_{\kappa,\mathcal{R}(p,q)}[y]_{n-\kappa,\mathcal{R}(p,q)}}\nonumber\\
			\nonumber\\
			&+& \sum_{\kappa=0}^{n-1}\bigg[\begin{array}{c} n-1  \\ \kappa\end{array} \bigg]_{\mathcal{R}(p,q)}\epsilon_1^{(\kappa+1)(y-n+\kappa+1)}\epsilon_2^{(n-\kappa-1)(x-\kappa)}{[x]_{\kappa+1,\mathcal{R}(p,q)}[y]_{n-\kappa-1,\mathcal{R}(p,q)}}\nonumber\\
			&=& \big[x+y-n+1\big]_{\mathcal{R}(p,q)}\,T_{n-1}(x;y)_{\mathcal{R}(p,q)}.
		\end{eqnarray*}
		Therefore, for $n\in \mathbb{N}\backslash \{0\},$ the sum $T_n(x;y)_{\mathcal{R}(p,q)}$ satisfies the first-order recursion relation 
		\begin{equation*}
		T_n(x;y)_{\mathcal{R}(p,q)} =  [x+y-n+1]_{\mathcal{R}(p,q)}\,T_{n-1}(x;y)_{\mathcal{R}(p,q)},
		\end{equation*}
		with $T_1(x;y)_{\mathcal{R}(p,q)} = [x+y]_{\mathcal{R}(p,q)}.$ Recursively,  
		it follows that $T_n(x;y)_{\mathcal{R}(p,q)} =
		{[x+y]_{n,\mathcal{R}(p,q)}}.$
		Therefore, we get  (\ref{vd1}). Finally, interchanging $x$ by $y,$ and  replacing $\kappa$ by $n-\kappa,$ the expression (\ref{vd1}) is rewritten in the form (\ref{vd2}). $\cqfd$
		\begin{remark}
			Note that the $q-$ deformed Vandermonde's formula can be retrieved by taking $\mathcal{R}(u,v)=v:$
			\begin{equation*}
			{[x+y]_{n,q}}
			= \sum_{\kappa=0}^{n}q^{(n-\kappa)(x-\kappa)}\bigg[{n \atop \kappa}\bigg]_q\,{[x]_{\kappa,q}\,[y]_{n-\kappa,q}},
			\end{equation*}	
			or, in an alternative form,
			\begin{equation*}
			{[x+y]_{n,q}}
			= \sum_{\kappa=0}^{n}q^{\kappa(y-n+\kappa)}\bigg[{n \atop \kappa}\bigg]_q\,{[x]_{\kappa,q}\,[y]_{n-\kappa,q}},
			\end{equation*}
			where $0<q<1.$
		\end{remark}
		From the $\mathcal{R}(p,q)-$ Vandermonde's formula (\ref{vd1}),  we can  deduce  the following remarkable $\mathcal{R}(p,q)-$ deformed identities.
		\begin{lemma}
			Let $x,$ $y,$ $p,$ and $q$ be real numbers such that $0<q<p\leq 1.$ Then, the following relations hold.
			\begin{equation}\label{ii1}
			{{[x+y+n]_{n,\mathcal{R}(p,q)}\over [y+n]_{n,\mathcal{R}(p,q)}}}
			= \sum_{\kappa=0}^{n}\bigg[\begin{array}{c} n  \\ \kappa\end{array} \bigg]_{\mathcal{R}(p,q)}\epsilon_1^{\kappa(y+\kappa)}\epsilon_2^{(n-\kappa)(x-\kappa)}{{[x]_{\kappa,\mathcal{R}(p,q)}\over [y+\kappa]_{\kappa,\mathcal{R}(p,q)}},}
			\end{equation}	
			\begin{equation}\label{ii2}
			{1\over \bigg[\begin{array}{c} x-1  \\ n \end{array} \bigg]_{\mathcal{R}(p,q)}}
			= \sum_{\kappa=0}^{n}(-1)^{n-\kappa}\bigg[\begin{array}{c} n  \\ \kappa\end{array} \bigg]_{\mathcal{R}(p,q)}\epsilon_1^{{n-\kappa\choose 2}+\kappa(n-x)}\,\epsilon_2^{{n-\kappa\choose 2}}{[x]_{\mathcal{R}(p,q)}\over [x-\kappa]_{\mathcal{R}(p,q)}}
			\end{equation}
			and
			\begin{equation}\label{ii3}
			{1\over \bigg[\begin{array}{c} y+n  \\ n \end{array} \bigg]^{-1}_{\mathcal{R}(p,q)}}= \sum_{\kappa=0}^{n}(-1)^\kappa\bigg[\begin{array}{c} n  \\ \kappa\end{array} \bigg]_{\mathcal{R}(p,q)}\epsilon_1^{{\kappa+1\choose 2}}\,\epsilon_2^{{\kappa+1\choose 2}-n(y+\kappa)}{[y]_{\mathcal{R}(p,q)}\over [y+\kappa]_{\mathcal{R}(p,q)}}.
			\end{equation}
		\end{lemma}
		{\it Proof:}
		Replacing $y$ by $y+n$ in (\ref{vd1}), we obtain
		\begin{equation*}
		{[x+y+n]_{n,\mathcal{R}(p,q)}}
		= \sum_{\kappa=0}^{n}\bigg[\begin{array}{c} n  \\ \kappa\end{array} \bigg]_{\mathcal{R}(p,q)}\epsilon_1^{\kappa(y+\kappa)}\epsilon_2^{(n-\kappa)(x-\kappa)}{[x]_{\kappa,\mathcal{R}(p,q)}[y+n]_{n-\kappa,\mathcal{R}(p,q)}.}
		\end{equation*}
		Multiplying both sides of this relation by
		{$[y]_{-n,\mathcal{R}(p,q)},$}
		and using
		{$[y]_{-\kappa,\mathcal{R}(p,q)}=[y]_{-n,\mathcal{R}(p,q)}[y+n]_{n-\kappa,\mathcal{R}(p,q)},$}
		we get
		\begin{equation*}
		{[y]_{-n,\mathcal{R}(p,q)}[x+y+n]_{n,\mathcal{R}(p,q)}}
		= \sum_{\kappa=0}^{n}\bigg[\begin{array}{c} n  \\ \kappa\end{array} \bigg]_{\mathcal{R}(p,q)}\epsilon_1^{\kappa(y+\kappa)}\epsilon_2^{(n-\kappa)(x-\kappa)}{[x]_{\kappa,\mathcal{R}(p,q)}\,[y]_{-\kappa,\mathcal{R}(p,q)}}
		\end{equation*}
		and according to (\ref{08}), we deduce the required formula
		\begin{equation*}
		{{[x+y+n]_{n,\mathcal{R}(p,q)}\over [y+n]_{n,\mathcal{R}(p,q)}}}
		= \sum_{\kappa=0}^{n}\bigg[\begin{array}{c} n  \\ \kappa\end{array} \bigg]_{\mathcal{R}(p,q)}\epsilon_1^{\kappa(y+\kappa)}\epsilon_2^{(n-\kappa)(x-\kappa)}{{[x]_{\kappa,\mathcal{R}(p,q)}\over [y+\kappa]_{\kappa,\mathcal{R}(p,q)}}.}
		\end{equation*}
		Putting now $y=-x$ in (\ref{ii1}), and using, respectively,
		\begin{equation*}
		{[-x+\kappa]_{\kappa,\mathcal{R}(p,q)}}
		=(-1)^\kappa(\epsilon_1\epsilon_2)^{{\kappa\choose 2}-\kappa\,x}
		{[x-1]_{\kappa,\mathcal{R}(p,q)}}
		\end{equation*}
		and
		\begin{equation*}
		{n-\kappa\choose 2}={n\choose 2}-{\kappa\choose 2}-\kappa(n-\kappa),
		\end{equation*}
		we get (\ref{ii2}). Similarly, by substituting $x$ by $-y,$ we obtain (\ref{ii3}). $\cqfd$
		\begin{remark}
			Taking $\mathcal{R}(p,q)=q,$ we  retrieve the $q-$ identities as particular cases:
			\begin{equation}\label{qii1}
			{{[u+v+n]_{n,q}\over [v+n]_{n,q}}}
			= \sum_{\kappa=0}^{n}\bigg[\begin{array}{c} n  \\ \kappa\end{array} \bigg]_{q}q^{(n-\kappa)(u-\kappa)}{{[u]_{\kappa,q}\over [v+\kappa]_{\kappa,q}},}
			\end{equation}	
			\begin{equation}\label{qii2}
			{1\over \bigg[\begin{array}{c} u-1  \\ n \end{array} \bigg]_{q}}= \sum_{\kappa=0}^{n}(-1)^{n-\kappa}\bigg[\begin{array}{c} n  \\ \kappa\end{array} \bigg]_{q}q^{{n-\kappa\choose 2}}{[u]_{q}\over [u-\kappa]_{q}}
			\end{equation}
			and
			\begin{equation}\label{qii3}
			{1\over\bigg[\begin{array}{c} v+n  \\ n \end{array} \bigg]_{q}}= \sum_{\kappa=0}^{n}(-1)^\kappa\bigg[\begin{array}{c} n  \\ \kappa\end{array} \bigg]_{q}q^{{\kappa+1\choose 2}-n(v+\kappa)}{[v]_{q}\over [v+\kappa]_{q}}.
			\end{equation}
		\end{remark}
		Considering two real numbers $u$ and  $v$  leads to  the following  results:
		\begin{theorem}
			The $\mathcal{R}(p,q)-$ deformed Cauchy's formula is given by: 
			\begin{equation}\label{c1}
			\bigg[\begin{array}{c} u+v  \\ n\end{array} \bigg]_{\mathcal{R}(p,q)} = \sum_{\kappa=0}^{n}\epsilon_1^{\kappa(v-n+\kappa)}\epsilon_2^{(n-\kappa)(u-\kappa)}\bigg[\begin{array}{c} u  \\ \kappa\end{array} \bigg]_{\mathcal{R}(p,q)}\bigg[\begin{array}{c} v  \\ n-\kappa \end{array} \bigg]_{\mathcal{R}(p,q)}
			\end{equation}
			or, equivalently,
			\begin{equation}\label{c2}
			\bigg[\begin{array}{c} u+v  \\ n\end{array} \bigg]_{\mathcal{R}(p,q)} = \sum_{\kappa=0}^{n}\epsilon_1^{(n-\kappa)(u-\kappa)}\epsilon_2^{\kappa(v-n+\kappa)}\bigg[\begin{array}{c} u  \\ \kappa\end{array} \bigg]_{\mathcal{R}(p,q)}\bigg[\begin{array}{c} v  \\ n-\kappa \end{array} \bigg]_{\mathcal{R}(p,q)}.
			\end{equation}
		\end{theorem}
		{\it Proof:} From (\ref{010}) and the $\mathcal{R}(p,q)-$ deformed Vandermonde's formula, we get the result.
		$\cqfd$
		
		We obtain the $q-$ deformed Cauchy's formulae (\ref{qc1}) and (\ref{qc2}) by taking $\mathcal
		{R}(p,q)=q.
		$
		\begin{theorem}
			The negative $\mathcal{R}(p,q)-$deformed Vandermonde's formula is given by:	\begin{equation}\label{nvd1}
			{[u+v]_{-n,\mathcal{R}(p,q)}}
			= \sum_{\kappa=0}^{\infty}\bigg[\begin{array}{c} -n  \\ \kappa\end{array} \bigg]_{\mathcal{R}(p,q)}\epsilon_1^{\kappa(v+n+\kappa)}\epsilon_2^{(-n-\kappa)(u-\kappa)}{[u]_{\kappa,\mathcal{R}(p,q)}[v]_{-n-\kappa,\mathcal{R}(p,q)}}
			\end{equation}
			or
			\begin{equation}\label{nvd2}
			{[u+v]_{-n,\mathcal{R}(p,q)}}
			= \sum_{\kappa=0}^{\infty}\bigg[\begin{array}{c} -n  \\ \kappa\end{array} \bigg]_{\mathcal{R}(p,q)}\epsilon_1^{(-n-\kappa)(u-\kappa)}\epsilon_2^{\kappa(v+n+\kappa)}{[u]_{\kappa,\mathcal{R}(p,q)}[v]_{-n-\kappa,\mathcal{R}(p,q)}}
			\end{equation}
			where 
			$n$ is a positive integer.
		\end{theorem}
		{\it Proof:} For $n\in\mathbb{N}\backslash\{0\},$ we consider the following  expression: 
		\begin{equation*}\label{nvd3}
		H_n(u;v)_{\mathcal{R}(p,q)}= \sum_{\kappa=0}^{\infty}\bigg[\begin{array}{c} -n  \\ \kappa\end{array} \bigg]_{\mathcal{R}(p,q)}\epsilon_1^{\kappa(v+n+\kappa)}\epsilon_2^{(-n-\kappa)(u-\kappa)}{[u]_{\kappa,\mathcal{R}(p,q)}[v]_{-n-\kappa,\mathcal{R}(p,q)}}.
		\end{equation*}
		For $n=1,$ we have
		\begin{equation*}\label{nvd6}
		H_1(u;v)_{\mathcal{R}(p,q)}= {1\over \big[u+v+1\big]_{\mathcal{R}(p,q)}}.
		\end{equation*}
		Using  the relation (\ref{bc1}) and  
		\begin{eqnarray*}\label{nvd8}
			\big[u+v+n+1\big]_{\mathcal{R}(p,q)}{[u]_{\kappa,\mathcal{R}(p,q)}[v]_{-n-\kappa-1,\mathcal{R}(p,q)}}
			&=& \epsilon_1^{v+n+\kappa+1}
			{[u]_{\kappa+1,\mathcal{R}(p,q)}[v]_{-n-\kappa-1,\mathcal{R}(p,q)}}
			\nonumber\\ &+& \epsilon_2^{u-\kappa}{[u]_{\kappa,\mathcal{R}(p,q)}[v]_{-n-\kappa,\mathcal{R}(p,q)}},
		\end{eqnarray*} 
		$H_n(u;v,p,q)$ takes the following form
		\begin{eqnarray*}\label{nvd7}
			H_n(u;v)_{\mathcal{R}(p,q)} &=& \sum_{\kappa=0}^{\infty}\bigg[\begin{array}{c} -n-1  \\ \kappa\end{array} \bigg]_{\mathcal{R}(p,q)}\epsilon_1^{\kappa(v+n+\kappa)+\kappa}\epsilon_2^{(-n-\kappa)(u-\kappa)}{[u]_{\kappa,\mathcal{R}(p,q)}[v]_{-n-\kappa,\mathcal{R}(p,q)}}\nonumber\\
			&+& \sum_{\kappa=0}^{\infty}\bigg[\begin{array}{c} -n-1  \\ \kappa\end{array} \bigg]_{\mathcal{R}(p,q)}\epsilon_1^{\kappa(v+n+\kappa)+\kappa}\epsilon_2^{(-n-\kappa)(u-\kappa+1)}{[u]_{\kappa+1,\mathcal{R}(p,q)}[v]_{-n-\kappa-1,\mathcal{R}(p,q)}}\nonumber\\
			&=& {H_{n-1}(u;v)_{\mathcal{R}(p,q)} \over [u+v+n+1]_{\mathcal{R}(p,q)}}.
		\end{eqnarray*}
		Therefore, for $n\in\mathbb{N}\backslash\{0\},$ the sum $H_n(u;v)_{\mathcal{R}(p,q)}$ satisfies the first-order recursion relation 
		\begin{equation*}\label{nvd10}
		H_n(u;v)_{\mathcal{R}(p,q)} ={H_{n-1}(u;v)_{\mathcal{R}(p,q)} \over \big[u+v+n+1\big]_{\mathcal{R}(p,q)}},\quad n\in\mathbb{N}\backslash\{0\},
		\end{equation*}
		with $H_1(u;v)_{\mathcal{R}(p,q)} =  {1\over [u+v+1]_{\mathcal{R}(p,q)}}.$ 
		Recursively, it comes that $H_n(u;v)_{p,q} = {[u+v]_{-n,\mathcal{R}(p,q)}}.$ 
		Following the steps used to prove (\ref{nvd1}), 
		we obtain (\ref{nvd2}), and the proof is achieved. $\cqfd$
		
		We recover the negative $q-$ Vandermonde's formulae  by taking $\mathcal{R}(p,q)=q$ as follows:
		\begin{equation*}
		{[x+y]_{-n,q}}
		= \sum_{\kappa=0}^{n}q^{-(n-\kappa)(x-\kappa)}\bigg[{-n \atop \kappa}\bigg]_q\,{[x]_{\kappa,q}\,[y]_{-n-\kappa,q}},
		\quad |q^{-(x+y+1)}|<1,
		\end{equation*}	
		and, alternatively, 
		\begin{equation*}
		{[x+y]_{-n,q}}
		= \sum_{\kappa=0}^{n}q^{\kappa)(y+n+\kappa)}\bigg[{-n \atop \kappa}\bigg]_q\,{[x]_{\kappa,q}\,[y]_{-n-\kappa,q}},
		\quad |q^{-(x+y+1)}|<1,
		\end{equation*}
		where $0<q< 1.$
		\begin{lemma}
			\begin{equation}\label{ndv11}
			{[v]_{-n,\mathcal{R}(p,q)}}
			= \sum_{\kappa=0}^{\infty}\bigg[\begin{array}{c} n+ \kappa -1  \\ \kappa\end{array} \bigg]_{\mathcal{R}(p,q)}\epsilon^{n(u-\kappa)}_1\epsilon^{\kappa(v-n+1)}_2{[u]_{\kappa,\mathcal{R}(p,q)}\over [u+v]_{n+\kappa,\mathcal{R}(p,q)}}
			\end{equation}
			and
			\begin{equation}\label{ndv12}
			{[v]_{-n,\mathcal{R}(p,q)}}
			= \sum_{\kappa=0}^{\infty}\bigg[\begin{array}{c} n+ \kappa -1  \\ \kappa\end{array} \bigg]_{\mathcal{R}(p,q)}\epsilon^{\kappa(v-n+1)}_1\epsilon^{n(u-\kappa)}_2{[u]_{\kappa,\mathcal{R}(p,q)}\over [u+v]_{n+\kappa,\mathcal{R}(p,q)}}.
			\end{equation}
		\end{lemma}
		{\it Proof:} For $n$ a positive integer, the $\mathcal{R}(p,q)-$ deformed factorial of $x=-n$ of order $\kappa$ is written as:
		\begin{equation*}\label{ndv13}
		{[-n]_{\kappa,\mathcal{R}(p,q)}}
		= (-1)^{\kappa}\,\Big(\epsilon_1\epsilon_2\Big)^{-n\kappa- {\kappa \choose 2}}\,{[n+\kappa-1]_{\kappa,\mathcal{R}(p,q)}}
		\end{equation*}
		and 
		\begin{equation*}\label{ndv14}
		\mathcal{R}(p^{-j},q^{-j}) =- \Big(\epsilon_1\,\epsilon_2\Big)^{-j}\,\mathcal{R}(p^j,q^j)\mbox{,}\quad j\in\{n,n+1,\cdots, n+\kappa-1\}.
		\end{equation*}
		In the same vein, the $\mathcal{R}(p,q)-$ deformed binomial  coefficient of $x=-n$ is given by:
		\begin{equation*}\label{ndv15}
		\bigg[\begin{array}{c} -n  \\ \kappa \end{array} \bigg]_{\mathcal{R}(p,q)} =(-1)^{\kappa}\Big(\epsilon_1\,\epsilon_2\Big)^{-n\kappa -{\kappa \choose 2} }\,\,\bigg[\begin{array}{c} n+\kappa -1  \\ \kappa \end{array} \bigg]_{\mathcal{R}(p,q)}.
		\end{equation*}
		Moreover,
		\begin{eqnarray*}\label{ndv16}
			{{ [-v-1]_{-n,\mathcal{R}(p,q)}}}
			&=& {1\over[-v-1+n]_{n,\mathcal{R}(p,q)}}
			\nonumber\\
			&=&{1 \over (-1)^n\,\Big(\epsilon_1\,\epsilon_2\Big)^{-nv + {n \choose 2}}{[v]_{n,\mathcal{R}(p,q)}}}
		\end{eqnarray*} 
		and 
		\begin{eqnarray*}\label{ndv17}
			{[-u-v-1]_{-n-\kappa,\mathcal{R}(p,q)}}
			&=& {1\over[-u-v-1+n+\kappa]_{n+\kappa,\mathcal{R}(p,q)}}
			\nonumber\\
			&=&{(-1)^{-n-\kappa} \over \Big(\epsilon_1\,\epsilon_2\Big)^{-(n+\kappa)(u+v) + {n+\kappa \choose 2}}{[u+v]_{n+\kappa,\mathcal{R}(p,q)}}}.
		\end{eqnarray*} 
		The relation (\ref{nvd1}) may be  written as follows:
		\begin{equation}\label{ndv18}
		{1\over (\epsilon_1\,\epsilon_2)^{-nv+{n\choose 2}}{[v]_{n,\mathcal{R}(p,q)}}}
		= \sum_{\kappa=0}^{\infty}\bigg[\begin{array}{c} n+\kappa -1  \\ \kappa \end{array} \bigg]_{\mathcal{R}(p,q)}\,\,F(p,q)
		\end{equation} 
		where
		\begin{equation*}\label{ndv19}
		F(p,q) = {\epsilon^{{\kappa(v+n+\kappa)}}_1(\epsilon_1\epsilon_2)^{-(n+\kappa)(u-\kappa) -n\kappa- {\kappa \choose 2}}{[u]_{\kappa,\mathcal{R}(p,q)}}
			\over (\epsilon_1\,\epsilon_2)^{-(n+\kappa)(u+v) + {n+\kappa \choose 2}}{[u+v]_{n+\kappa,\mathcal{R}(p,q)}}}.
		\end{equation*}
		Taking into acount the fact that
		\begin{equation*}\label{ndv20}
		{n +\kappa \choose 2} =  {n  \choose 2} + {\kappa \choose 2} + n\kappa,
		\end{equation*}
		the relation (\ref{ndv18}) is reduced to (\ref{ndv11}). Similarly, we obtain (\ref{ndv12}). 
		$\cqfd$
		
		Setting $\epsilon_1=1$ and $\epsilon_2=q$ provides the $q-$analogs of the formulae (\ref{ndv11}) and (\ref{ndv12}) as:
		\begin{equation*}\label{qndv11}
		{[v]_{-n,q}}
		= \sum_{\kappa=0}^{\infty}\bigg[\begin{array}{c} n+ \kappa -1  \\ \kappa\end{array} \bigg]_{q}q^{\kappa(v-n+1)}{{[u]_{\kappa,q}\over [u+v]_{n+\kappa,q}}},
		\quad |q^{v}|<1
		\end{equation*}
		and
		\begin{equation*}\label{qndv12}
		{[v]_{-n,q}}
		= \sum_{\kappa=0}^{\infty}\bigg[\begin{array}{c} n+ \kappa -1  \\ \kappa\end{array} \bigg]_{q}q^{n(u-\kappa)}{{[u]_{\kappa,q}\over [u+v]_{n+\kappa,q}}},
		\quad |q^{-v}|<1.
		\end{equation*}
		\begin{remark}
			By taking $\mathcal{R}(u,v)=((q-p^{-1})v)^{-1}(uv-1),$ 
			and $n-$  a positive integer, we deduct, as particular cases,
			the \textbf{generalized  $q-$Quesne} deformed Vandermonde and Cauchy's formulae given, respectively,  by:
			\begin{eqnarray*}\label{vdga}
				{[u+v]^Q_{n,p,q}}
				&=& \sum_{\kappa=0}^{n}\bigg[\begin{array}{c} n  \\ \kappa\end{array} \bigg]^Q_{p,q}p^{\kappa(v-n+\kappa)}q^{-(n-\kappa)(u-\kappa)}{[u]^Q_{\kappa,p,q}\,[v]^Q_{n-\kappa,p,q}}
				\\
				&=& \sum_{\kappa=0}^{n}\bigg[\begin{array}{c} n  \\ \kappa\end{array} \bigg]^Q_{p,q}p^{(n-\kappa)(u-\kappa)}q^{-\kappa(v-n+\kappa)}{[u]^Q_{\kappa,p,q}\,[v]^Q_{n-\kappa,p,q}}
			\end{eqnarray*}
			and
			\begin{eqnarray*}
				\bigg[\begin{array}{c} u+v  \\ n\end{array} \bigg]^Q_{p,q} &=& \sum_{\kappa=0}^{n-1}p^{\kappa(v-n+\kappa)}q^{-(n-\kappa)(u-\kappa)}\bigg[\begin{array}{c} u  \\ \kappa\end{array} \bigg]^Q_{p,q}\bigg[\begin{array}{c} v  \\ n-\kappa \end{array} \bigg]^Q_{p,q}\\
				& =& \sum_{\kappa=0}^{n-1}p^{(n-\kappa)(u-\kappa)}q^{-\kappa(v-n+\kappa)}\bigg[\begin{array}{c} u  \\ \kappa\end{array} \bigg]^Q_{p,q}\bigg[\begin{array}{c} v  \\ n-\kappa \end{array} \bigg]^Q_{p,q},
			\end{eqnarray*}
			while their negative counterparts
			are
			provided by
			\begin{eqnarray*}\label{nvd1g}
				{[u+v]^Q_{-n,p,q}}
				&=&
				\sum_{\kappa=0}^{\infty}\bigg[\begin{array}{c} -n  \\ \kappa\end{array} \bigg]^Q_{p,q}p^{\kappa(v+n+\kappa)}q^{-(-n-\kappa)(u-\kappa)}{[u]^Q_{\kappa,p,q}\,[v]^Q_{-n-\kappa,p,q}}
				\\
				&=& \sum_{\kappa=0}^{\infty}\bigg[\begin{array}{c} -n  \\ \kappa\end{array} \bigg]^Q_{p,q}p^{(-n-\kappa)(u-\kappa)}q^{-\kappa(v+n+\kappa)}{[u]^Q_{\kappa,p,q}\,[v]^Q_{-n-\kappa,p,q}}
			\end{eqnarray*}
			and
			\begin{eqnarray*}\label{ndv11g}
				{[v]^Q_{-n,p,q}}
				&=& \sum_{\kappa=0}^{\infty}\bigg[\begin{array}{c} n+ \kappa -1  \\ \kappa\end{array} \bigg]^Q_{p,q}\,p^{n(u-\kappa)}\,q^{-\kappa(v-n+1)}{{[u]^Q_{\kappa,p,q}\over [u+v]^Q_{n+\kappa,p,q}}}
				\\
				& = &\sum_{\kappa=0}^{\infty}\bigg[\begin{array}{c} n+ \kappa -1  \\ \kappa\end{array} \bigg]^Q_{p,q}\,p^{\kappa(v-n+1)}\,q^{-n(u-\kappa)}{{[u]^Q_{\kappa,p,q}\over [u+v]^Q_{n+\kappa,p,q}}},
			\end{eqnarray*}
			respectively.
		\end{remark}
		\subsection{$\mathcal{R}(p,q)-$ deformed binomial and negative  binomial formulae}
		In this section, we examine in detail the $\mathcal{R}(p,q)-$ deformed binomial. Here also,
		$x,$ $p,$ and $q$ are  real numbers with $ 0 < q < p \leq 1;$ $n$ is a positive integer. Then, 
		\begin{theorem}
			\begin{equation}\label{bn1}
			\prod_{r=1}^{n}\Big(\epsilon_1^{r-1} + x\,\epsilon_2^{r-1}\Big) = \sum_{\kappa=0}^{n}\epsilon_1^{{n-\kappa \choose 2}}\,\epsilon_2^{{\kappa \choose 2}}\, \bigg[\begin{array}{c} n  \\ \kappa \end{array} \bigg]_{\mathcal{R}(p,q)}\,x^{\kappa} 
			\end{equation}
		\end{theorem}
		{\it Proof:} The result follows from  induction on $n.$  $\cqfd$
		
		Taking $\mathcal{R}(p,q)=q,$ we recover the $q-$ binomial formula (\ref{qbc1}),
		while  $\mathcal{R}(u,v)=(p-q)^{-1}(u-v)$  gives the $(p,q)-$ binomial formula (\ref{l1}).
		\begin{theorem}
			\begin{equation}\label{bn5}
			\prod_{r=1}^{n}\Big(\epsilon_1^{r-1} - x\,\epsilon_2^{r-1}\Big)^{-1} = \epsilon^{-{n \choose 2}}_1\sum_{\kappa=0}^{\infty}\epsilon^{-\kappa(n-1)}_1 \bigg[\begin{array}{c} n + \kappa -1  \\ \kappa \end{array} \bigg]_{\mathcal{R}(p,q)}\,x^{\kappa}
			\end{equation}
			where 
			$|\epsilon_2|<|\epsilon_1|.$
		\end{theorem}
		{\it Proof:} We have
		\begin{equation*}
		\prod_{r=1}^{n}\Big(\epsilon_1^{r-1} - x\,\epsilon_2^{r-1}\Big)^{-1}=\epsilon^{-{n \choose 2}}_1\prod_{r=1}^{n}\Big(1 - x\big({\epsilon_2\over \epsilon_1}\big)^{r-1}\Big)^{-1}.
		\end{equation*} 
		Setting $\theta={\epsilon_2\over \epsilon_1}<1$ and from \cite{CA1}, we obtain
		\begin{eqnarray*}
			\prod_{r=1}^{n}\Big(\epsilon_1^{r-1} - x\,\epsilon_2^{r-1}\Big)^{-1}&=&\epsilon^{-{n \choose 2}}_1\sum_{\kappa=0}^{\infty} \bigg[\begin{array}{c} n + \kappa -1  \\ \kappa \end{array} \bigg]_{\mathcal{R}({q\over p})}\,x^{\kappa}\nonumber\\
			&=& \epsilon^{-{n \choose 2}}_1\sum_{\kappa=0}^{\infty}\epsilon^{-\kappa(n-1)}_1 \bigg[\begin{array}{c} n + \kappa -1  \\ \kappa \end{array} \bigg]_{\mathcal{R}(p,q)}\,x^{\kappa}.
		\end{eqnarray*}
		The proof is achieved. $\cqfd$
		
		The negative $q-$ binomial coefficient can be obtained, by setting $\mathcal{R}(u,v)=v,$ in the form:
		\begin{equation*}\label{qbn}
		\prod_{r=1}^{n}\Big(1 - t\,q^{r-1}\Big)^{-1} =\sum_{\kappa=0}^{\infty} \bigg[\begin{array}{c} n + \kappa -1  \\ \kappa \end{array} \bigg]_{q}\,t^{\kappa}.
		\end{equation*}
		A novel  negative $\mathcal{R}(p,q)-$ binomial formula can be deduced as follows:
		\begin{lemma}
			For  $0<x<\infty,$ 
			\begin{equation}\label{ad1}
			\prod_{i=1}^{n}(\epsilon^{i-1}_1+x\epsilon^{i-1}_2)=\epsilon^{n\choose 2}_1\sum_{\kappa=0}^{\infty}\bigg[\begin{array}{c} n + \kappa -1  \\ \kappa \end{array} \bigg]_{\mathcal{R}(p,q)}\,{x^\kappa\,\epsilon^\kappa_1\,\epsilon^{\kappa\choose 2}_2\over \displaystyle\prod_{i=1}^{\kappa}(\epsilon^{n+i-1}_1+x\epsilon^{n+i-1}_2)}
			\end{equation}
			or, in an equivalent way,
			\begin{equation}\label{ad2}
			{\displaystyle\prod_{i=1}^{n}(\epsilon^{i-1}_1+x\epsilon^{i-1}_2)\over x^n\,\epsilon^{n\choose 2}_2}=\sum_{\kappa=0}^{\infty}\bigg[\begin{array}{c} n + \kappa -1  \\ \kappa \end{array} \bigg]_{\mathcal{R}(p,q)}\,{\epsilon^\kappa_2\,\epsilon^{\kappa\choose 2}_1\over \displaystyle\prod_{i=1}^{\kappa}(\epsilon^{n+i-1}_1+x\epsilon^{n+i-1}_2)}.
			\end{equation}
		\end{lemma}
		{\it Proof:} Since
		\begin{equation*}
		\prod_{i=1}^{n}\Big(\epsilon_1^{i-1} - x\,\epsilon_2^{i-1}\Big)=\epsilon^{{n \choose 2}}_1\prod_{i=1}^{n}\Big(1 - x\big({\epsilon_2\over \epsilon_1}\big)^{i-1}\Big)^{-1},
		\end{equation*} 
		and using \cite{CA1}, we obtain
		\begin{eqnarray*}
			\prod_{i=1}^{n}\Big(\epsilon_1^{i-1} - x\,\epsilon_2^{i-1}\Big)&=&\epsilon^{{n \choose 2}}_1\sum_{\kappa=0}^{\infty} \bigg[\begin{array}{c} n + \kappa -1  \\ \kappa \end{array} \bigg]_{\mathcal{R}({q\over p})}\,{x^{\kappa}\,({\epsilon_2\over \epsilon_1})^{\kappa\choose 2}\over \displaystyle\prod_{i=1}^{\kappa}\Big(1+x({\epsilon_2\over \epsilon_1})^{n+i-1}\Big)}.
		\end{eqnarray*}
		From (\ref{a016}) and 
		\begin{equation*}
		\displaystyle\prod_{i=1}^{\kappa}\Big(1+x({\epsilon_2\over \epsilon_1})^{n+i-1}\Big)=\epsilon^{{-\kappa(n-1)}-{\kappa+1\choose 2}}_1\,\,\displaystyle\prod_{i=1}^{\kappa}\Big(\epsilon^{n+i-1}_1+x\,\epsilon^{n+i-1}_2\Big)
		\end{equation*}
		we get (\ref{ad1}). Replacing $x$ by $x^{-1},$ $\epsilon_1$ by $\epsilon^{-1}_1$ and $\epsilon_2$ by $\epsilon^{-1}_2$ leads to (\ref{ad2}),  and the proof is achieved. $\cqfd$
		\begin{remark}
			The  negative $q-$binomial formula can be obtained by taking $\mathcal{R}(u,v)=v:$
			\begin{equation}\label{qad1}
			\prod_{i=1}^{n}(1+x\,q^{i-1})=\sum_{\kappa=0}^{\infty}\bigg[\begin{array}{c} n + \kappa -1  \\ \kappa \end{array} \bigg]_{q}\,{x^\kappa\,q^{\kappa\choose 2}\over \displaystyle\prod_{i=1}^{\kappa}(1+x\,q^{n+i-1})}
			\end{equation}
			which can  also  be translated into the form:
			\begin{equation}\label{qad2}
			{\displaystyle\prod_{i=1}^{n}(1+x\,q^{i-1})\over x^n\,q^{n\choose 2}}=\sum_{\kappa=0}^{\infty}
			\bigg[\begin{array}{c} n + \kappa -1  \\ \kappa \end{array} \bigg]_{q}
			\,{q^\kappa \over \displaystyle\prod_{i=1}^{\kappa}(1+xq^{n+i-1})}.
			\end{equation}
		\end{remark} 
		\begin{theorem}
			The following orthogonality relations hold:
			\begin{equation}\label{bn6}
			\sum_{x=\kappa}^{n}(-1)^{n-x}\epsilon_1^{{x \choose 2}}\,\epsilon_2^{{n-x \choose 2}}\bigg[\begin{array}{c} n   \\ x \end{array} \bigg]_{\mathcal(p,q)}\,\bigg[\begin{array}{c} x   \\ \kappa \end{array} \bigg]_{\mathcal{R}(p,q)}=\delta_{n,\kappa},
			\end{equation}
			and
			\begin{equation}\label{bn7}
			\sum_{x=\kappa}^{n}(-1)^{x-\kappa}\epsilon_1^{{\kappa \choose 2}}\,\epsilon_2^{{x-\kappa \choose 2}}\bigg[\begin{array}{c} n   \\ x \end{array} \bigg]_{\mathcal{R}(p,q)}\,\bigg[\begin{array}{c} x   \\ \kappa \end{array} \bigg]_{\mathcal{R}(p,q)}=\delta_{n,\kappa},
			\end{equation}
			where  $\delta_{n,\kappa}$ is the Kronecker delta, $n$ and $\kappa$ are positive integers.
		\end{theorem}
		{\it Proof:} Since 
		\begin{eqnarray*}\label{bn8}
			\bigg[\begin{array}{c} n   \\ x \end{array} \bigg]_{\mathcal{R}(p,q)}\,\bigg[\begin{array}{c} x   \\ \kappa \end{array} \bigg]_{\mathcal{R}(p,q)}
			&=& \bigg[\begin{array}{c} n   \\ \kappa \end{array} \bigg]_{\mathcal{R}(p,q)}\,\,\bigg[\begin{array}{c} n-\kappa   \\ n-x \end{array} \bigg]_{\mathcal{R}(p,q)},
		\end{eqnarray*}
		then, the relation (\ref{bn6}) may be expressed as:
		\begin{eqnarray*}\label{bn9}
			\sum_{x=\kappa}^{n}(-1)^{n-x}\epsilon_1^{{x \choose 2}}\,\epsilon_2^{{n-x \choose 2}}\bigg[\begin{array}{c} n   \\ x \end{array} \bigg]_{\mathcal{R}(p,q)}\,\bigg[\begin{array}{c} x   \\ \kappa \end{array} \bigg]_{\mathcal{R}(p,q)}
			&=& \delta_{n,\kappa}.
		\end{eqnarray*}
		From the expression
		\begin{equation*}\label{bn10}
		\bigg[\begin{array}{c} n   \\ x \end{array} \bigg]_{\mathcal{R}(p,q)}\,\bigg[\begin{array}{c} x   \\ \kappa \end{array} \bigg]_{\mathcal{R}(p,q)}= \bigg[\begin{array}{c} n   \\ \kappa \end{array} \bigg]_{\mathcal{R}(p,q)}\,\,\bigg[\begin{array}{c} n-\kappa   \\ x-\kappa \end{array} \bigg]_{\mathcal{R}(p,q)},
		\end{equation*}
		we get (\ref{bn7}). Therefore, the result holds. 
		$\cqfd$
		
		The $q-$ deformed orthogonality relations (\ref{qbn6}) and (\ref{qbn7}) can be recovered by putting $\mathcal{R}(p,q)=q.$
		\begin{corollary}
			The  inversion of the $\mathcal{R}(p,q)-$ binomial formula is provided by:
			\begin{equation}\label{bn11}
			\sum_{\kappa=0}^{n}(-1)^{\kappa}\,\epsilon_1^{{\kappa \choose 2}}\,\epsilon_2^{{\kappa \choose 2}} \bigg[\begin{array}{c} n   \\ \kappa \end{array} \bigg]_{\mathcal{R}(p,q)}\prod_{r=1}^{\kappa}\Big(\epsilon_1^{1-r} - x\,\epsilon_2^{1-r}\Big)=x^n.
			\end{equation}
			In particular,
			\begin{equation}\label{bna11}
			\sum_{\kappa=0}^{n}(-1)^{\kappa}\,\epsilon_1^{{\kappa \choose 2}-x\kappa}\,\epsilon_2^{{\kappa \choose 2}} \bigg[\begin{array}{c} n   \\ \kappa \end{array} \bigg]_{\mathcal{R}(p,q)}(\epsilon_1-\epsilon_2)^\kappa\,{[x]_{\kappa,\mathcal{R}(p,q)}}= \big({\epsilon_2 \over \epsilon_1}\big)^{nx}.
			\end{equation}
		\end{corollary}
		{\it Proof:} From the relation (\ref{016}) and replacing $\epsilon_1$ by $\epsilon_1^{-1},$ $\epsilon_2$ by $\epsilon_2^{-1},$ $x$ by $-x$, $n$ by $\kappa$ and 
		$\kappa$ by $r$ in (\ref{bn1}), we get 
		\begin{eqnarray*}\label{bn12}
			\prod_{i=1}^{\kappa}\Big(\epsilon_1^{1-i} - x\,\epsilon_2^{1-i}\Big)
			&=& \sum_{r=0}^{\kappa}(-1)^r\epsilon_1^{-{\kappa-r\choose 2}-r(\kappa-r)}\epsilon_2^{-{r \choose 2}-r(\kappa-r)}\,\bigg[\begin{array}{c} n   \\ \kappa \end{array} \bigg]_{\mathcal{R}(p,q)}x^r.
		\end{eqnarray*}
		Multiplying the members of the above expression by $(-1)^{\kappa}\,\epsilon_1^{{\kappa \choose 2}}\,\epsilon_2^{{\kappa \choose 2}}\bigg[\begin{array}{c} n   \\ \kappa \end{array} \bigg]_{\mathcal{R}(p,q)},$ it comes
		\begin{eqnarray}\label{bn13}
		C_1:&=&
		(-1)^{\kappa}\epsilon_1^{{\kappa \choose 2}}\epsilon_2^{{\kappa \choose 2}}\bigg[\begin{array}{c} n   \\ \kappa \end{array} \bigg]_{\mathcal{R}(p,q)}\prod_{i=1}^{\kappa}\Big(\epsilon_1^{1-i} - x\epsilon_2^{1-i}\Big)
		\nonumber\\
		&=&\sum_{r=0}^{\kappa}\bigg[\begin{array}{c} n   \\ \kappa \end{array} \bigg]_{\mathcal{R}(p,q)}(-1)^{\kappa-r}\epsilon_1^{-{\kappa-r\choose 2}-r(\kappa-r)+{\kappa\choose 2}}
		\epsilon_2^{-{r \choose 2}-r(\kappa-r)+{\kappa\choose 2}}\bigg[\begin{array}{c} n   \\ \kappa \end{array} \bigg]_{\mathcal{R}(p,q)}x^r.
		\end{eqnarray}
		Summing the expression (\ref{bn13}) for $\kappa\in\{1,2,\cdots,n\},$  using (\ref{bn7}) and
		\begin{equation}\label{bn14}
		{\kappa-r\choose 2} = {\kappa\choose 2} - {r\choose 2} -r(\kappa-r),
		\end{equation}
		we obtain (\ref{bn11}). Replacing $x$ by $\epsilon^{-x}_1\,\epsilon^x_2$ in (\ref{bn11}) and using 
		\begin{equation*}\label{bn15}
		\prod_{r=1}^{\kappa}\Big(\epsilon^{1-r}_1 - \epsilon^{-x}_1\,\epsilon^{x+1-r}_2\Big)=\epsilon^{-\kappa\,x}_1 \big(\epsilon_1-\epsilon_2\big)^\kappa{[x]_{\kappa,\mathcal{R}(p,q)}}
		\end{equation*} 
		yield (\ref{bna11}).  $\cqfd$
		\begin{remark}
			Putting $\epsilon_1=1$ and $\epsilon_2=q,$ we obtain the  inversion of the $q-$ binomial formula \cite{CA1}:
			\begin{equation*}\label{qbn11}
			\sum_{\kappa=0}^{n}(-1)^{\kappa}\,q^{{\kappa \choose 2}} \bigg[\begin{array}{c} n   \\ \kappa \end{array} \bigg]_{q}\prod_{r=1}^{\kappa}(1 - t\,q^{1-r})=t^n.
			\end{equation*}
			In particular,
			\begin{equation*}\label{qbna11}
			\sum_{\kappa=0}^{n}(-1)^{\kappa}\,q^{{\kappa \choose 2}} \bigg[\begin{array}{c} n   \\ \kappa \end{array} \bigg]_{q}(1-q)^\kappa\,{[t]_{\kappa,q}}= q^{nt}.
			\end{equation*}
		\end{remark}
		The inversion of the deformed binomial formulae (\ref{bn11}) and (\ref{bna11}) leads to the following results:
		\begin{lemma}
			\begin{equation}\label{bn16}
			\sum_{\kappa=0}^{n} \bigg[\begin{array}{c} n   \\ \kappa \end{array} \bigg]_{\mathcal{R}(p,q)}\,x^{n-\kappa}\prod_{r=1}^{\kappa}\Big(\epsilon_1^{r-1} - x\,\epsilon_2^{r-1}\Big)=1
			\end{equation}
			and
			\begin{equation*}\label{bna17}
			\sum_{\kappa=0}^{n} \bigg[\begin{array}{c} n   \\ \kappa \end{array} \bigg]_{\mathcal{R}(p,q)}\,\epsilon_1^{-\kappa(n-\kappa)}\,\epsilon_2^{-\kappa(x+n -\kappa) }(\epsilon_1-\epsilon_2)^\kappa\,{[x]_{\kappa,\mathcal{R}(p,q)}}=\big({\epsilon_1\over \epsilon_2}\big)^{nx}.
			\end{equation*}
		\end{lemma}
		{\it Proof:} Replacing $x$ by $x^{-1}$ in equation (\ref{bn11}),
		and using the expression
		\begin{equation*}\label{bn18}
		(-1)^\kappa\,\epsilon_1^{{\kappa \choose 2}}\,\epsilon_2^{{\kappa \choose 2}}\prod_{r=1}^{\kappa}\Big(\epsilon_1^{1-r} - x^{-1}\,\epsilon_2^{1-r}\Big)=x^{-\kappa}\prod_{r=1}^{\kappa}\Big(\epsilon_1^{r-1} - x\,\epsilon_2^{r-1}\Big),
		\end{equation*}
		we obtain (\ref{bn16}). Also replacing $\epsilon_1$ by $\epsilon^{-1}_1$ and $\epsilon_2$ by $\epsilon^{-1}_2$ in (\ref{bna11}),  we get
		\begin{equation*}
		\sum_{\kappa=0}^{n}(-1)^{\kappa}\,\epsilon_1^{-{\kappa \choose 2}+ x\kappa}\,\epsilon_2^{-{\kappa \choose 2}} \bigg[\begin{array}{c} n   \\ \kappa \end{array} \bigg]_{\mathcal{R}(p^{-1},q^{-1})}(\epsilon^{-1}_1-\epsilon^{-1}_2)^\kappa\,{[x]_{\kappa,\mathcal{R}(p^{-1},q^{-1})}}= \big({\epsilon_1 \over \epsilon_2}\big)^{nx}.
		\end{equation*}
		Using the relations (\ref{014}), (\ref{016}) 
		and after computation, the result holds. $\cqfd$
		
		It is worth noticing the following relevant   identities from  $\mathcal{R}(p,q)-$ binomial and  negative $\mathcal{R}(p,q)-$ binomial formulae, as exposed and proved below:
		\begin{itemize}
			\item
			\begin{lemma}
				Let $r,$ $s$ and $n$ be positive integers. Then, the following relations hold:
				\begin{equation}\label{a11}
				\sum_{\kappa=0}^{n}\epsilon^{\kappa(s-n+\kappa)}_1\epsilon^{(n-\kappa)(r-\kappa)}_2\bigg[\begin{array}{c} r   \\ \kappa \end{array} \bigg]_{\mathcal{R}(p,q)}\bigg[\begin{array}{c} s   \\ n-\kappa \end{array} \bigg]_{\mathcal{R}(p,q)}=\bigg[\begin{array}{c} r+s   \\ n \end{array} \bigg]_{\mathcal{R}(p,q)},
				\end{equation}
				\begin{equation}\label{a12}
				\sum_{\kappa=0}^{n}\epsilon^{\kappa(m+\kappa)}_1\epsilon^{(r-\kappa)(r-\kappa-m)}_2\bigg[\begin{array}{c} r   \\ \kappa \end{array} \bigg]_{\mathcal{R}(p,q)}\bigg[\begin{array}{c} r   \\ \kappa+m \end{array} \bigg]_{\mathcal{R}(p,q)}=\bigg[\begin{array}{c} 2r   \\ \kappa+m \end{array} \bigg]_{\mathcal{R}(p,q)}
				\end{equation}
				and
				\begin{equation}\label{a13}
				\sum_{\kappa=0}^{r}\epsilon^{\kappa^2}_1\epsilon^{(r-\kappa)^2}_2\bigg[\begin{array}{c} r   \\ \kappa \end{array} \bigg]^2_{\mathcal{R}(p,q)}
				=\bigg[\begin{array}{c} 2r   \\ r \end{array} \bigg]_{\mathcal{R}(p,q)}.
				\end{equation}	
			\end{lemma}
			{\it Proof:} From the $\mathcal{R}(p,q)-$ binomial formula and 
			\begin{equation*}
			\prod_{i=1}^{r}(\epsilon^{i-1}_1+x\epsilon^{i-1}_2)\prod_{i=1}^{s}(\epsilon^{r+i-1}_1+x\epsilon^{r+i-1}_2)=\prod_{i=1}^{r+s}(\epsilon^{i-1}_1+x\epsilon^{i-1}_2),
			\end{equation*}
			we get 
			\begin{equation*}
			\sum_{n=0}^{r+s}a_n\,x^n=\sum_{n=0}^{r+s}\epsilon^{r+s-n\choose 2}_1\epsilon^{n\choose 2}_2\bigg[\begin{array}{c} r+s   \\ n \end{array} \bigg]_{\mathcal{R}(p,q)}x^n,
			\end{equation*}
			where
			\begin{equation*}
			a_n\,x^n = \sum_{\kappa=0}^{r}\epsilon^{r-\kappa\choose 2}_1\epsilon^{\kappa \choose 2}_2\bigg[\begin{array}{c} r   \\ \kappa \end{array} \bigg]_{\mathcal{R}(p,q)}x^\kappa\sum_{\kappa=0}^{s}\epsilon^{{s-j\choose 2}+r(s-j)}_1\epsilon^{{j \choose 2}+r\,j}_2\bigg[\begin{array}{c} s   \\ j \end{array} \bigg]_{\mathcal{R}(p,q)}x^j.
			\end{equation*}
			Using
			\begin{equation*}
			{n-\kappa\choose 2}={n\choose 2}-{\kappa\choose 2}-\kappa(n-\kappa)
			\end{equation*}
			and
			\begin{equation*}
			{r+s-n\choose 2}={r-\kappa \choose 2}+{s-n+\kappa\choose 2}+(r-\kappa)(s-n+\kappa),
			\end{equation*}
			we obtain
			\begin{equation*}
			a_n\,x^n=\epsilon^{r+s-n\choose 2}_1\epsilon^{n\choose 2}_2\sum_{\kappa=0}^{n}\epsilon^{\kappa(s-n+\kappa)}_1\epsilon^{(n-\kappa)(r-\kappa)}_2\bigg[\begin{array}{c} r   \\ \kappa \end{array} \bigg]_{\mathcal{R}(p,q)}\bigg[\begin{array}{c} s   \\ n-\kappa \end{array} \bigg]_{\mathcal{R}(p,q)}
			\end{equation*}
			giving (\ref{a11}). Setting $s=r$ and $n=r-m$ in (\ref{a11}), and using
			\begin{equation*}
			\bigg[\begin{array}{c} r   \\ r-m-\kappa \end{array} \bigg]_{\mathcal{R}(p,q)}=\bigg[\begin{array}{c} r   \\ \kappa +m \end{array} \bigg]_{\mathcal{R}(p,q)},\mbox{and}\,\bigg[\begin{array}{c} 2r   \\ r-m \end{array} \bigg]_{\mathcal{R}(p,q)}=\bigg[\begin{array}{c} 2r   \\ r+m \end{array} \bigg]_{\mathcal{R}(p,q)}
			\end{equation*}
			the relation (\ref{a12}) holds. Since
			\begin{equation*}
			\bigg[\begin{array}{c} r   \\ r-\kappa \end{array} \bigg]_{\mathcal{R}(p,q)}=\bigg[\begin{array}{c} r   \\ \kappa \end{array} \bigg]_{\mathcal{R}(p,q)},
			\end{equation*}
			and putting $r=n=s$ in (\ref{a11}), we get (\ref{a13}), what achieves the proof. $\cqfd$
			\begin{remark} Note the following:
				\begin{enumerate}
					\item[(1)] The relation (\ref{a11}) is a particular case of the $\mathcal{R}(p,q)-$deformed Cauchy formula (\ref{c1}), with $u=r$
					and $v=s.$
					\item[(2)] The particular case of  $q-$ deformation  is achieved  from the above formulae by taking $\mathcal{R}(p,q)=q$ as follows:
					\begin{equation*}\label{qa11}
					\sum_{\kappa=0}^{n}q^{(n-\kappa)(r-\kappa)}\bigg[\begin{array}{c} r   \\ \kappa \end{array} \bigg]_{q}\bigg[\begin{array}{c} s   \\ n-\kappa \end{array} \bigg]_{q}=\bigg[\begin{array}{c} r+s   \\ n \end{array} \bigg]_{q},
					\end{equation*}
					\begin{equation*}\label{qa12}
					\sum_{\kappa=0}^{n}q^{(r-\kappa)(r-\kappa-m)}\bigg[\begin{array}{c} r   \\ \kappa \end{array} \bigg]_{q}\bigg[\begin{array}{c} r   \\ \kappa+m \end{array} \bigg]_{q}=\bigg[\begin{array}{c} 2r   \\ \kappa+m \end{array} \bigg]_{q}
					\end{equation*}
					and
					\begin{equation*}\label{qa13}
					\sum_{\kappa=0}^{r}q^{(r-\kappa)^2}\bigg[\begin{array}{c} r   \\ \kappa \end{array} \bigg]^2_{q}
					=\bigg[\begin{array}{c} 2r   \\ r \end{array} \bigg]_{q}.
					\end{equation*}
				\end{enumerate}
			\end{remark}
			\item
			\begin{lemma}\label{L327}
				Let $r,$ $s,$ and $n$ be positive integers. Then, 
				\begin{eqnarray*}\label{a21}
					\bigg[\begin{array}{c} r+s+n-1   \\ n \end{array} \bigg]_{\mathcal{R}(p,q)}=\sum_{\kappa=0}^{n}\epsilon^{r(n-s-\kappa)+\kappa\,s}_1\epsilon^{r(n-\kappa)}_2\bigg[\begin{array}{c} r+\kappa-1   \\ \kappa \end{array} \bigg]_{\mathcal{R}(p,q)}
					\bigg[\begin{array}{c} s+n-\kappa-1   \\ n-\kappa \end{array} \bigg]_{\mathcal{R}(p,q)}.
				\end{eqnarray*}
			\end{lemma}
			{\it Proof:}
			From the relation
			\begin{equation*}
			\prod_{i=1}^{r}(\epsilon^{i-1}_1+x\epsilon^{i-1}_2)^{-1}\prod_{i=1}^{s}(\epsilon^{r+i-1}_1+x\epsilon^{r+i-1}_2)^{-1}=\prod_{i=1}^{r+s}(\epsilon^{i-1}_1+x\epsilon^{i-1}_2)^{-1},
			\end{equation*}
			where $|x|<1,$ and from the negative $\mathcal{R}(p,q)-$binomial formula, we get
			\begin{equation*}
			\sum_{n=0}^{\infty}b_n\,x^n=\sum_{n=0}^{\infty}\epsilon^{-{r+s \choose 2}-n(r+s-1)}_1\bigg[\begin{array}{c} r+s+n-1   \\ n \end{array} \bigg]_{\mathcal{R}(p,q)}x^n,
			\end{equation*}
			where
			\begin{eqnarray*}
				&&b_n\,x^n = \sum_{\kappa=0}^{\infty}\epsilon^{{r\choose 2}-\kappa(r-1)}_1\bigg[\begin{array}{c} r+\kappa-1   \\ \kappa \end{array} \bigg]_{\mathcal{R}(p,q)}x^\kappa\cr&&\qquad\qquad\qquad\qquad\qquad\times\sum_{j=0}^{\infty}\epsilon^{-{s\choose 2}-j(s-1)}_1\epsilon^{r\,j}_2\bigg[\begin{array}{c} s+j-1   \\ j \end{array} \bigg]_{\mathcal{R}(p,q)}x^j.
			\end{eqnarray*}
			Setting $j=n-\kappa,$ and after computation, 
			it follows
			\begin{eqnarray*}
				&&b_n\,x^n=\epsilon^{-{r+s \choose 2}-n(r+s-1)}_1\sum_{\kappa=0}^{n}\epsilon^{r(n-s-\kappa)+\kappa\,s}_1\epsilon^{r(n-\kappa)}_2\cr&&\qquad\qquad\qquad\qquad\qquad\times\bigg[\begin{array}{c} r+\kappa-1   \\ \kappa \end{array} \bigg]_{\mathcal{R}(p,q)}\bigg[\begin{array}{c} s+ n-\kappa -1  \\ n-\kappa \end{array} \bigg]_{\mathcal{R}(p,q)}
			\end{eqnarray*}
			yielding the result. $\cqfd$
			\begin{remark}
				\begin{enumerate}
					\item[(a)] This formula constitutes a particular version of the $\mathcal{R}(p,q)-$ Cauchy formula (\ref{c1}), corresponding to negative integers  $x=-r$ and $y=-s.$ 
					\item[(b)] Putting $\mathcal{R}(x,y)=y,$ we recover the $q-$ deformed analog of this Lemma as:
					\begin{eqnarray*}\label{qa21}
						\bigg[\begin{array}{c} r+s+n-1   \\ n \end{array} \bigg]_{q}=\sum_{\kappa=0}^{n}q^{r(n-\kappa)}\bigg[\begin{array}{c} r+\kappa-1   \\ \kappa \end{array} \bigg]_{q}
						\bigg[\begin{array}{c} s+n-\kappa-1   \\ n-\kappa \end{array} \bigg]_{q}.
					\end{eqnarray*}
				\end{enumerate}
			\end{remark}
			\item
			\begin{lemma}
				For $n,$ $m,$ and $\kappa,$  positive integers, such that $\kappa\leq m\leq n,$ we have:
				\begin{equation*}\label{a31}
				\bigg[\begin{array}{c} n   \\ m \end{array} \bigg]_{\mathcal{R}(p,q)}=\sum_{\kappa=0}^{n}\epsilon^{{(\kappa-2r-m+n)\kappa}+r(m+1)}_1\epsilon^{\kappa(n-r)}_2\bigg[\begin{array}{c} r-1   \\ \kappa-1 \end{array} \bigg]_{\mathcal{R}(p,q)}
				\bigg[\begin{array}{c} n-r   \\ m-\kappa \end{array} \bigg]_{\mathcal{R}(p,q)}.
				\end{equation*}
			\end{lemma}
			{\it Proof:}
			Multiplying both the members of the negative $\mathcal{R}(p,q)-$ binomial formula 
			\begin{equation*}
			\prod_{i=1}^{m+1}(\epsilon_1^{i-1}-x\epsilon_2^{i-1})^{-1}=\sum_{j=0}^{\infty}\epsilon^{-{m+1\choose 2}-m\,j}_1\bigg[\begin{array}{c} m+j   \\ m \end{array} \bigg]_{\mathcal{R}(p,q)}x^j
			\end{equation*}
			by $x^m$ and putting $n=m+j,$ we obtain
			\begin{equation*}
			x^m\prod_{i=1}^{m+1}(\epsilon_1^{i-1}-x\epsilon_2^{i-1})^{-1}=\sum_{n=m}^{\infty}\epsilon_1^{-{m+1\choose 2}-m(n-m)}\bigg[\begin{array}{c} n   \\ m \end{array} \bigg]_{\mathcal{R}(p,q)}x^n.
			\end{equation*}
			Similarly, we get
			\begin{equation*}
			x^\kappa\prod_{i=1}^{\kappa}(\epsilon_1^{i-1}-x\epsilon_2^{i-1})^{-1}=\sum_{r=\kappa}^{\infty}\epsilon_1^{-{\kappa\choose 2}-(\kappa-1)(r-\kappa)}\bigg[\begin{array}{c} r-1   \\ \kappa-1 \end{array} \bigg]_{\mathcal{R}(p,q)}x^r
			\end{equation*}
			and 
			\begin{eqnarray*}
				x^{m-\kappa}\prod_{i=1}^{m-\kappa+1}(\epsilon_1^{\kappa+i-1}-x\epsilon_2^{\kappa+i-1})^{-1}&=&\sum_{j=m-\kappa}^{\infty}\epsilon_1^{-{m-\kappa+1\choose 2}-(m-\kappa)(j-m-\kappa)}\nonumber\\&\times&\epsilon^{\kappa\,j}_2\bigg[\begin{array}{c} j   \\ m-\kappa \end{array} \bigg]_{\mathcal{R}(p,q)}x^{m-\kappa}.
			\end{eqnarray*}
			Using
			\begin{equation*}
			x^{\kappa}\prod_{i=1}^{\kappa}(\epsilon^{\kappa+i-1}_1+x\epsilon^{i-1}_2)^{-1}x^{m-\kappa}\prod_{i=1}^{m-\kappa+1}(\epsilon^{\kappa+i-1}_1+x\epsilon^{\kappa+i-1}_2)^{-1}=x^m\prod_{i=1}^{m+1}(\epsilon^{i-1}_1+x\epsilon^{i-1}_2)^{-1}
			\end{equation*}
			and the above relations, we arrive at the result. $\cqfd$
			
			Note that the $q-$ deformed version can be recovered by taking $\mathcal{R}(p,q)=q:$
			\begin{equation*}\label{qa31}
			\sum_{\kappa=0}^{n}q^{\kappa(n-r)}\bigg[\begin{array}{c} r-1   \\ \kappa-1 \end{array} \bigg]_{q}
			\bigg[\begin{array}{c} n-r   \\ m-\kappa \end{array} \bigg]_{q}=\bigg[\begin{array}{c} n   \\ m \end{array} \bigg]_{q}.
			\end{equation*}
			\begin{remark} It is worth noticing that
				\begin{enumerate} 
					\item[(1)] The \textbf{generalized $q-$Quesne} deformed binomial and negative binomial formulae can straightforwardly  be retrieved from the general formalism as particular cases as follows: 
					\begin{itemize}
						\item[(i)] 
						\begin{equation*}\label{bn1g}
						\prod_{r=1}^{n}\Big(p^{r-1} + x\,q^{-r+1}\Big) = \sum_{\kappa=0}^{n}p^{{n-\kappa \choose 2}}\,q^{-{\kappa \choose 2}}\, \bigg[\begin{array}{c} n  \\ \kappa \end{array} \bigg]^Q_{p,q}\,x^{\kappa}. 
						\end{equation*}
						\item[(ii)]
						\begin{equation*}\label{bn5q}
						\prod_{r=1}^{n}\Big(p^{r-1} - x\,q^{1-r}\Big)^{-1} = p^{-{n \choose 2}}\sum_{\kappa=0}^{\infty} \bigg[\begin{array}{c} n + \kappa -1  \\ \kappa \end{array} \bigg]^Q_{p,q}\,x^{\kappa}.
						\end{equation*}
						
						\item[(iii)]
						\begin{equation*}\label{bn6g}
						\sum_{x=\kappa}^{n}(-1)^{n-x}p^{{x \choose 2}}\,q^{-{n-x \choose 2}}\bigg[\begin{array}{c} n   \\ x \end{array} \bigg]^Q_{p,q}\,\bigg[\begin{array}{c} x   \\ \kappa \end{array} \bigg]^Q_{p,q}=\delta_{n,\kappa},
						\end{equation*}
						and
						\begin{equation*}\label{bn7g}
						\sum_{x=\kappa}^{n}(-1)^{x-\kappa}p^{{\kappa \choose 2}}\,q^{-{x-\kappa \choose 2}}\bigg[\begin{array}{c} n   \\ x \end{array} \bigg]^Q_{p,q}\,\bigg[\begin{array}{c} x   \\ \kappa \end{array} \bigg]^Q_{p,q}=\delta_{n,\kappa}.
						\end{equation*}
						\item[(iv)]
						\begin{equation}\label{bn11g}
						\sum_{\kappa=0}^{n}(-1)^{\kappa}\,p^{{\kappa \choose 2}}\,q^{-{\kappa \choose 2}} \bigg[\begin{array}{c} n   \\ \kappa \end{array} \bigg]^Q_{p,q}\prod_{r=1}^{\kappa}\Big(p^{1-r} - x\,q^{-1+r}\Big)=x^n.
						\end{equation}
						In particular,
						\begin{equation}\label{bna11g}
						\sum_{\kappa=0}^{n}(-1)^{\kappa}\,p^{{\kappa \choose 2}-x\kappa}\,q^{-{\kappa \choose 2}} \bigg[\begin{array}{c} n   \\ \kappa \end{array} \bigg]^Q_{p,q}(q-p^{-1})^\kappa\,{[x]^Q_{\kappa,p,q}}= (pq)^{-nx}.
						\end{equation}
					\end{itemize}
					\item[(2)] The associated
					deformed binomial formulae (\ref{bn11g}) and (\ref{bna11g}) can also be rewritten as follows:
					\begin{equation*}\label{bn21}
					\sum_{\kappa=0}^{n} \bigg[\begin{array}{c} n   \\ \kappa \end{array} \bigg]^Q_{p,q}\,x^{n-\kappa}\prod_{r=1}^{\kappa}\Big(p^{r-1} - x\,q^{-r+1}\Big)=1
					\end{equation*}
					and
					\begin{equation*}\label{bna22}
					\sum_{\kappa=0}^{n} \bigg[\begin{array}{c} n   \\ \kappa \end{array} \bigg]^Q_{p,q}\,p^{-\kappa(n-\kappa)}\,q^{\kappa(x+n-\kappa) }(q-p^{-1})^\kappa\,{[x]^Q_{\kappa,p,q}}=\big(p\,q\big)^{nx}.
					\end{equation*}
					Further, the following identities are retrieved from the above Lemmas:
					\begin{equation*}\label{ga11}
					\sum_{\kappa=0}^{n}p^{\kappa(s-n+\kappa)}\,q^{-(n-\kappa)(r-\kappa)}\bigg[\begin{array}{c} r   \\ \kappa \end{array} \bigg]^Q_{p,q}\bigg[\begin{array}{c} s   \\ n-\kappa \end{array} \bigg]^Q_{p,q}=\bigg[\begin{array}{c} r+s   \\ n \end{array} \bigg]^Q_{p,q},
					\end{equation*}
					\begin{equation*}\label{ga12}
					\sum_{\kappa=0}^{n}p^{\kappa(m+\kappa)}\,q^{-(r-\kappa)(r-\kappa-m)}\bigg[\begin{array}{c} r   \\ \kappa \end{array} \bigg]^Q_{p,q}\bigg[\begin{array}{c} r   \\ \kappa+m \end{array} \bigg]^Q_{p,q}=\bigg[\begin{array}{c} 2r   \\ \kappa+m \end{array} \bigg]^Q_{p,q},
					\end{equation*}
					\begin{equation*}\label{ga13}
					\sum_{\kappa=0}^{r}p^{\kappa^2}\,q^{-(r-\kappa)^2}\bigg(\bigg[\begin{array}{c} r   \\ \kappa \end{array} \bigg]^Q_{p,q}\bigg)^2
					=\bigg[\begin{array}{c} 2r   \\ r \end{array} \bigg]^Q_{p,q},
					\end{equation*}
					\begin{eqnarray*}\label{ga21}
						\bigg[\begin{array}{c} r+s+n-1   \\ n \end{array} \bigg]^Q_{p,q}&=&\sum_{\kappa=0}^{n}p^{r(n-s-\kappa)+\kappa\,s}\,q^{-r(n-\kappa)}\bigg[\begin{array}{c} r+\kappa-1   \\ \kappa \end{array} \bigg]^Q_{p,q}\nonumber\\&\times&
						\bigg[\begin{array}{c} s+n-\kappa-1   \\ n-\kappa \end{array} \bigg]^Q_{p,q}
					\end{eqnarray*}
					\begin{equation*}\label{ga31}
					\bigg[\begin{array}{c} n   \\ m \end{array} \bigg]^Q_{p,q}=\sum_{\kappa=0}^{n}p^{{(\kappa-2r-m+n)\kappa}+r(m+1)}\,q^{-\kappa(n-r)}\bigg[\begin{array}{c} r-1   \\ \kappa-1 \end{array} \bigg]^Q_{p,q}
					\bigg[\begin{array}{c} n-r   \\ m-\kappa \end{array} \bigg]^Q_{p,q}
					\end{equation*}
					\begin{equation*}\label{gad1}
					\prod_{i=1}^{n}(p^{i-1}+x\,q^{-i+1})=p^{n\choose 2}\sum_{\kappa=0}^{\infty}\bigg[\begin{array}{c} n + \kappa -1  \\ \kappa \end{array} \bigg]^Q_{p,q}\,{x^\kappa\,p^\kappa\,q^{-{\kappa\choose 2}} \over \displaystyle\prod_{i=1}^{\kappa}(p^{n+i-1}+x\,q^{-n-i+1})}
					\end{equation*}
					or
					\begin{equation*}\label{gad2}
					{\displaystyle\prod_{i=1}^{n}(p^{i-1}+x\,q^{-i+1})\over x^n\,q^{-{n\choose 2}}}=\sum_{\kappa=0}^{\infty}\bigg[\begin{array}{c} n + \kappa -1  \\ \kappa \end{array} \bigg]^Q_{p,q}\,{q^{-\kappa}\,p^{\kappa\choose 2}\over \displaystyle\prod_{i=1}^{\kappa}(p^{n+i-1}+x\,q^{-n-i+1})}.
					\end{equation*}
				\end{enumerate}
			\end{remark}
		\end{itemize}
		\subsection{$\mathcal{R}(p,q)-$ deformed Stirling numbers}
		\begin{definition}
			Let $x,$ $j,$  $p,$ and $q$ be  real numbers such that $0<q<p\leq 1.$ Then, the  noncentral $\mathcal{R}(p,q)-$ deformed factorial of $x$ of order $n$ and of noncentrality parameter $j$ is defined by:
			\begin{equation}\label{sn1}
			{[x-j]_{n,\mathcal{R}(p,q)}}
			:= [x-j]_{\mathcal{R}(p,q)}\,[x-j-1]_{\mathcal{R}(p,q)},\,\cdots,[x-j-n+1]_{\mathcal{R}(p,q)},
			\end{equation}
			where $n$ is a positive integer.
		\end{definition}
		Taking $j=0$ in (\ref{sn1}),  we obtain the $\mathcal{R}(p,q)-$ deformed factorial of $x$ of order $n.$
		Following the relation
		\begin{equation*}\label{s2}
		\mathcal{R}(p^{x-j-r},q^{x-j-r}) = \epsilon_2^{j+r}\Big(\mathcal{R}(p^x,q^x) - \epsilon_1^{x-j-r}\mathcal{R}(p^{j+r},q^{j+r})\Big)\mbox{,}\quad r\in\mathbb{N}\backslash\{0\},
		\end{equation*}
		the equation (\ref{sn1}) takes the  form
		\begin{eqnarray*}\label{s3}
			{[x-j]_{n,\mathcal{R}(p,q)}}
			&=&\epsilon_2^{-{n\choose 2}-j\,n}\Big([x]_{\mathcal{R}(p,q)} - \epsilon_1^{x-j}\,[j]_\mathcal{R}(p,q)\Big)\nonumber\\&\cdots& \Big([x]_{\mathcal{R}(p,q)} - \epsilon_1^{x-j-n+1}\,[j+n-1]_{\mathcal{R}(p,q)}\Big).
		\end{eqnarray*}
		Further, we get  a polynomial of the $\mathcal{R}(p,q)-$ deformed number of degree $n$ as follows:
		\begin{equation}\label{s4}
		{[x-j]_{n,\mathcal{R}(p,q)}}
		=\epsilon_2^{-{n\choose 2}-j\,n}\,\sum_{\kappa=0}^{n}s_{\mathcal{R}(p,q)}(n,\kappa;j)\,{[x]^\kappa_{\mathcal{R}(p,q)}}
		\mbox{,}\quad n\in\mathbb{N}
		\end{equation}
		or 
		\begin{equation}\label{s5}
		{[x]^n_{\mathcal{R}(p,q)}}
		=\sum_{\kappa=0}^{n}\epsilon_2^{{\kappa\choose 2}+j\,\kappa}\,S_{\mathcal{R}(p,q)}(n,\kappa;j)\,{[x-j]_{\kappa,\mathcal{R}(p,q)}}
		\mbox{,}\quad n\in\mathbb{N}
		\end{equation}
		Equivalently,
		\begin{equation}\label{s6}
		{[x+j]^n_{\mathcal{R}(p,q)}}
		=\sum_{\kappa=0}^{n}\epsilon_2^{{\kappa\choose 2}+j\,\kappa}\,S_{\mathcal{R}(p,q)}(n,\kappa;j)\,{[x]_{\kappa,\mathcal{R}(p,q)}}
		\mbox{,}\quad n\in\mathbb{N}.
		\end{equation}
		The coefficients $s_{\mathcal{R}(p,q)}(n, \kappa;j)$ and $S_{\mathcal{R}(p,q)}(n,\kappa;j)$ are  the noncentral $\mathcal{R}(p,q)-$ deformed Stirling numbers of
		the first and second kind, respectively.
		\begin{remark}
			\begin{enumerate}
				\item[(1)] 
				For $j = 0,$ these deformed numbers are
				reduced to $s_{\mathcal{R}(p,q)}(n, \kappa;0) = s_{\mathcal{R}(p,q)}(n, \kappa)$ and $S_{\mathcal{R}(p,q)}(n, \kappa; 0) =
				S_{\mathcal{R}(p,q)}(n, \kappa),$ which are nothing but the $\mathcal{R}(p,q)-$ deformed Stirling numbers of the first and second kinds, respectively.
				\item[(2)] Taking $\mathcal{R}(x,y)=y,$ we recover the $q-$ deformed Stirling numbers of first and second kinds   (\ref{qsnf}) and (\ref{qsns}) .
			\end{enumerate}
		\end{remark}
		\begin{lemma}
			\begin{equation*}\label{s7}
			|s_{\mathcal{R}(p^{-1},q^{-1})}(n, \kappa,j)| = \Bigg(\mathcal{R}(p^{-1},q^{-1})\Bigg)^{n-\kappa}\,\,s_{\mathcal{R}(p^{-1},q^{-1})}(n, \kappa,j)
			\end{equation*}
			and 
			\begin{equation*}\label{s8}
			|s_{\mathcal{R}(p,q)}(n, \kappa,j)| = \Bigg(\mathcal{R}(p^{-1},q^{-1})\Bigg)^{\kappa-n}\,\,s_{\mathcal{R}(p,q)}(n, \kappa,j),
			\end{equation*}
			where $n$ is a positive interger, and $|s_{\mathcal{R}(p,q)}(n, \kappa,j)|$ stands for  the deformed absolute noncentral $\mathcal{R}(p,q)-$ deformed Stirling number of the
			first kind.
		\end{lemma}
		{\it Proof:} Since 
		\begin{equation*}\label{s9}
		\mathcal{R}(p^{x+j},q^{x+j}) = \epsilon_2^{j}\Big(\mathcal{R}(p^x,q^x) + \epsilon_1^{x+j}\,(\epsilon_1\,\epsilon_2)^{-1}\mathcal{R}(p^j,q^j)\Big),
		\end{equation*}
		then,
		\begin{eqnarray*}\label{s10}
			[x+j+n-1]_{\mathcal{R}(p,q)}&=&\epsilon_2^{{n\choose 2}+jn}\Big([x]_{\mathcal{R}(p,q)}+\epsilon_1^{x+j-1}\epsilon_2^{-1}[j]_{\mathcal{R}(p,q)}\Big)\nonumber\\&\cdots&\Big([x]_{\mathcal{R}(p,q)} + \epsilon_1^{x+j+n-2}\epsilon_2^{-1}[j+n-1]_{\mathcal{R}(p,q)}\Big)
		\end{eqnarray*}
		is also a polynomial of the $\mathcal{R}(p,q)-$deformed number $[x]_{\mathcal{R}(p,q)}$ of degree $n,$ and after computation, for $n\in\mathbb{N}\backslash\{0\},$ we get
		\begin{equation*}\label{s11}
		{[x+j+n-1]_{n,\mathcal{R}(p,q)}}
		=\epsilon_2^{{n\choose 2}+jn}\sum_{\kappa=0}^{n}|s_{\mathcal{R}(p^{-1},q^{-1})}(n,\kappa,j)|\,\,
		[x]_{\mathcal{R}(p,q)}^{\kappa}.
		\end{equation*}
		Furthermore,
		\begin{equation*}\label{s12}
		{[x+j+n-1]_{n,\mathcal{R}(p,q)}}
		=\Bigg( \mathcal{R}(p^{-1},q^{-1})\Bigg)^n\,\,{[-x-j]_{n,\mathcal{R}(p^{-1},q^{-1})}}
		\end{equation*}
		\begin{equation*}
		\mathcal{R}(p^{-x},q^{-x}) =\mathcal{R}(p^{-1},q^{-1})\,\mathcal{R}(p^x,q^x).
		\end{equation*}
		Replacing now $-x$ by $x$ and $\epsilon_1^{-1}$ by $\epsilon_1,$ and $\epsilon_2^{-1}$ by $\epsilon_2$ in  (\ref{s4}), the result follows. $\cqfd$ 
		\begin{theorem}
			The $\mathcal{R}(p,q)-$ deformed Stirling numbers of the first and second kinds verify the orthogonality relations:
			\begin{equation}\label{s13}
			\sum_{m=\kappa}^{n}s_{\mathcal{R}(p,q)}(n,m,j)\,S_{\mathcal{R}(p,q)}(m,\kappa,j)= \delta_{n,\kappa}
			\end{equation}
			and 
			\begin{equation}\label{s14}
			\sum_{m=\kappa}^{n}S_{\mathcal{R}(p,q)}(n,m,j)\,s_{\mathcal{R}(p,q)}(m,\kappa,j)= \delta_{n,\kappa}.
			\end{equation}
		\end{theorem}
		{\it Proof:}  From the relations (\ref{s4}) and (\ref{s5}), we write
		\begin{eqnarray*}\label{s15}
			{[x-j]_{n,\mathcal{R}(p,q)}}
			&=&\epsilon_2^{-{n\choose 2}-jn}\sum_{m=0}^{n}s_{\mathcal{R}(p,q)}(n,\kappa;j){[x]^m_{\mathcal{R}(p,q)}}\nonumber\\ 
			&=&  \sum_{\kappa=0}^{n}\epsilon_2^{-{n\choose 2}+{\kappa\choose 2}-j(n-\kappa)}\nonumber\\&\times&\bigg(\sum_{m=\kappa}^{n}s_{p,q}(n,m;j)S_{\mathcal{R}(p,q)}(m,\kappa;j)\bigg){[x-j]_{\kappa,\mathcal{R}(p,q)}},
		\end{eqnarray*}
		giving, after computation,  (\ref{s13}). Similarly, we obtain (\ref{s14}). $\cqfd$
		
		The next statement is also valid.
		\begin{theorem}
			For $n\in\mathbb{N}\backslash\{0\}$ and $\kappa\in\{1,2,\cdots,n+1\},$ 
			the noncentral $\mathcal{R}(p,q)-$ deformed Stirling numbers of the first  and second kinds, $s_{\mathcal{R}(p,q)}$ and $S_{\mathcal{R}(p,q)},$ obey, respectively, the recursion relations
			\begin{eqnarray}\label{s16}
			s_{\mathcal{R}(p,q)}(n+1,\kappa;j)&=& s_{\mathcal{R}(p,q)}(n,\kappa-1;j)\nonumber\\
			&-& \epsilon^{x-n-j}_1[n+j]_{\mathcal{R}(p,q)}s_{\mathcal{R}(p,q)}(n,\kappa,j),\quad 
			\end{eqnarray}
			with initial conditions
			$s_{\mathcal{R}(p,q)}(0,0,j)=1,$  $s_{\mathcal{R}(p,q)}(n,0,j)=\epsilon_2^{{n \choose 2}+jn}\,{[-j]_{n,\mathcal{R}(p,q)}}
			\geq 0$ and $s_{\mathcal{R}(p,q)}(0,\kappa,j)=0,$ $j\geq 0;$ and
			%
			\begin{equation}\label{s17}
			S_{\mathcal{R}(p,q)}(n+1,\kappa;j)= S_{\mathcal{R}(p,q)}(n,\kappa-1;j)+ \epsilon^{x-\kappa}_1[\kappa+j]_{\mathcal{R}(p,q)}S_{\mathcal{R}(p,q)}(n,\kappa;j)
			\end{equation}
			with initial conditions $S_{\mathcal{R}(p,q)}(0,0,j)=1,$  $S_{\mathcal{R}(p,q)}(n,0,j)=[j]^n_{\mathcal{R}(p,q)},
			$
			$n\geq 0,$  and  $S_{\mathcal{R}(p,q)}(0,\kappa,j)=0,$ $\kappa\geq 0.$ 
		\end{theorem} 
		{\it Proof:} \begin{enumerate}
			\item[(1)] Let us consider the relation
			\begin{equation}\label{s18}
			{[x-j]_{n+1,\mathcal{R}(p,q)}}
			= \epsilon_2^{-n-j}\Big([x]_{\mathcal{R}(p,q)}-\epsilon^{x-n-j}_1[n+j]_{\mathcal{R}(p,q)}\Big)\,{[x-j]_{n,\mathcal{R}(p,q)}}
			\end{equation}
			or 
			\begin{equation*}
			{[x-j]_{n+1,\mathcal{R}(p,q)}}
			= \epsilon_1^{-n-j}\Big([x]_{\mathcal{R}(p,q)}-\epsilon^{x-n-j}_2\,[n+j]_{\mathcal{R}(p,q)}\Big)\,{[x-j]_{n,\mathcal{R}(p,q)}}.
			\end{equation*}
			From the relation (\ref{s4}) and the expansion of both members of the recursion relation (\ref{s18}) into powers of  $\mathcal{R}(p^x,q^x),$ we get
			\begin{eqnarray*}\label{s19}
				S_1:&=&
				\epsilon^{-{n+1 \choose 2}-j(n+1)}_2\sum_{\kappa=0}^{n+1}s_{\mathcal{R}(p,q)(n+1,\kappa;j)}\,\,
				[x]^\kappa_{\mathcal{R}(p,q)}
				\nonumber\\&=&\epsilon^{-{n \choose 2}-j(n+1)-n}_2\sum_{r=0}^{n}s_{\mathcal{R}(p,q)(n,r;j)}\,\,
				[x]^{r+1}_{\mathcal{R}(p,q)}
				\nonumber\\ &-&\epsilon^{-{n \choose 2}-j(n+1)-n}_2\,\epsilon^{x-n-j}_1\sum_{\kappa=0}^{n}s_{\mathcal{R}(p,q)(n,\kappa;j)}\,\,
				[x]^\kappa_{\mathcal{R}(p,q)}
				[n+j]_{\mathcal{R}(p,q)}\\
				&=&\epsilon^{-{n+1 \choose 2}-j(n+1)}_2\sum_{\kappa=1}^{n+1}s_{\mathcal{R}(p,q)(n,\kappa-1;j)}\,\,
				[x]^\kappa_{\mathcal{R}(p,q)}
				\nonumber\\ &-&\epsilon^{-{n+1 \choose 2}-j(n+1)}_2 \epsilon^{x-n-j}_1\sum_{\kappa=0}^{n}s_{\mathcal{R}(p,q)(n,\kappa;j)}\,\,
				[x]^\kappa_{\mathcal{R}(p,q)}
				[n+j]_{\mathcal{R}(p,q)}
			\end{eqnarray*}
			which gives (\ref{s16}). We use the relation (\ref{s4}) to get the initial conditions.
			\item[(2)] Similarly, consider
			\begin{eqnarray}\label{s21}
			\Bigg(\mathcal{R}(p^{x+j},q^{x+j})\Bigg)^{n+1}
			&=&	\mathcal{R}(p^{x+j},q^{x+j})	\Bigg(\mathcal{R}(p^{x+j},q^{x+j})\Bigg)^{n}
			\end{eqnarray} 
			and use the relation (\ref{s5}) to obtain
			\begin{eqnarray*}\label{s22}
				S_2:&=&\sum_{\kappa=0}^{n+1}\epsilon^{{\kappa \choose 2}+j\kappa}_2S_{\mathcal{R}(p,q)(n+1,\kappa;j)}{[x]_{\kappa,\mathcal{R}(p,q)}}
				\nonumber\\
				&=& \sum_{\kappa=0}^{n}\epsilon^{{\kappa \choose 2}+\kappa+j(\kappa+1)}_2\,S_{\mathcal{R}(p,q)(n,\kappa;j)}{[x]_{\kappa+1,\mathcal{R}(p,q)}}
				\nonumber\\
				&+& \sum_{\kappa=0}^{n}\epsilon^{{\kappa \choose 2}+j\,\kappa}_2\,\epsilon^{x-\kappa}_1S_{\mathcal{R}(p,q)(n,\kappa;j)}{[x]_{\kappa,\mathcal{R}(p,q)}}
				[\kappa+j]_{\mathcal{R}(p,q)}\\
				&=& \sum_{\kappa=1}^{n+1}\epsilon^{{\kappa \choose 2}+j\kappa}_2\,S_{\mathcal{R}(p,q)(n,\kappa-1;j)}{[x]_{\kappa,\mathcal{R}(p,q)}}
				\nonumber\\
				&+& \sum_{\kappa=0}^{n}\epsilon^{{\kappa \choose 2}+j\kappa}_2\,\epsilon^{x-\kappa}_1S_{\mathcal{R}(p,q)(n,\kappa;j)}{[x]_{\kappa,\mathcal{R}(p,q)}[\kappa+j]_{\mathcal{R}(p,q)}}
			\end{eqnarray*}
			Thus, we get the relation (\ref{s17}). The initial conditions follow from (\ref{s6}).  $\cqfd$
		\end{enumerate}
		\begin{theorem}
			For fixed $\kappa$, the generating function of  noncentral $\mathcal{R}(p,q)-$ deformed Stirling numbers of the second kind is given as follows:
			\begin{equation}\label{s24}
			\Psi_{\kappa}(v;p,q,r) = \sum_{n=\kappa}^{\infty}S_{\mathcal{R}(p,q)}(n;\kappa,r)\,v^n\mbox{,}\quad \kappa\in\mathbb{N}
			\end{equation}
			or, equivalently, in product form
			\begin{equation}\label{s25}
			\Psi_{\kappa}(v;p,q,r) = v^{\kappa}\,\prod_{j=0}^{\kappa}\Big(1 - \epsilon^{x-j}_1[r+j]_{\mathcal{R}(p,q)}\,v\Big)^{-1}
			\end{equation}
			for $|v|< {\epsilon^{-x+j}_1 \over [r+j]_{\mathcal{R}(p,q)}}.$
		\end{theorem}
		{\it Proof:} We assume that the series (\ref{s24}) converges. Multiplying the expression  (\ref{s17}) by $v^{n+1}$ and 
		summing the resulting relation for $n=\kappa-1,\kappa,\cdots,$ and $\kappa\in\mathbb{N}\backslash\{0\},$  we obtain
		\begin{equation*}
		\Psi_{\kappa}(v;p,q,r) =v\,\Psi_{\kappa-1}(v;p,q,r) + \epsilon^{x-\kappa}_1[r+\kappa]_{\mathcal{R}(p,q)}\,v\,\Psi_{\kappa}(v;p,q,r)\mbox{,}\quad 
		\end{equation*}
		which implies
		\begin{equation*}
		\Psi_{\kappa}(v;p,q,r) = v\,\Big( 1 -\epsilon^{x-\kappa}_1[r+\kappa]_{\mathcal{R}(p,q)}\,v \Big)^{-1}\,\Psi_{\kappa-1}(v;p,q,r).
		\end{equation*}
		By induction and consideration that 
		\begin{equation*}
		\Psi_{0}(v;p,q,r) = \sum_{n=0}^{\infty}S_{\mathcal{R}(p,q)}(n;0,r)\,v^n
		= \Big( 1 -\epsilon^{x}_1[r]_{\mathcal{R}(p,q)}\,v \Big)^{-1}
		\end{equation*}
		we obtain (\ref{s24}) and (\ref{s25}).$\cqfd$	
		\begin{remark}
			The generating function of the noncentral $q-$ Stirling numbers of the second kind can be obtained as:
			\begin{equation*}\label{qs24}
			\Psi_{\kappa}(v;q,r) = \sum_{n=\kappa}^{\infty}S_{q}(n;\kappa,r)\,v^n
			= v^{\kappa}\,\prod_{j=0}^{\kappa}(1 - [r+j]_{q}\,v)^{-1}
			\end{equation*}
			where $\kappa\in\mathbb{N}$ and $|v|< {1 \over [r+j]_{q}}.$
		\end{remark}
		\begin{lemma}
			For $\kappa\in\mathbb{N}$ and $n=\kappa, \kappa+1, \cdots, $  the reciprocal noncentral $\mathcal{R}(p,q)-$ deformed factorial ${
				{1 \over  [t-x]_{\kappa+1,\mathcal{R}(p,q)}}}$
			is expanded into the  reciprocal $\mathcal{R}(p,q)-$ deformed powers ${1\over [t]^{n+1}_{\mathcal{R}(p,q)}}$ as follows:
			\begin{equation}\label{s31}
			{\epsilon^{(t-r-x)\kappa}_1\over{
					[t-x]_{\kappa+1,\mathcal{R}(p,q)}} 
			}
			= \epsilon_2^{{\kappa +1\choose 2}+x(\kappa+1)}\sum_{n=\kappa}^{\infty}S_{\mathcal{R}(p,q)}(n,\kappa;x){\epsilon^{(t-r-x)n}_1\over [t]^{n+1}_{\mathcal{R}(p,q)}},
			\end{equation} 
			while
			the reciprocal $\mathcal{R}(p,q)-$ deformed powers  $[t]^{-\kappa-1}_{\mathcal{R}(p,q)}$  is expanded into   the reciprocal noncentral $\mathcal{R}(p,q)-$ deformed factorial ${{1\over [t-x]_{n+1,\mathcal{R}(p,q)}}}$ 
			as   below expressed:
			\begin{equation*}\label{s32}
			{\epsilon^{(t-r-x)\kappa}_1\over [t]^{\kappa+1}_{\mathcal{R}(p,q)}} =\sum_{n=\kappa}^{\infty} \epsilon_2^{-{n +1\choose 2}-x(n+1)}s_{\mathcal{R}(p,q)}(n,\kappa;x){\epsilon^{(t-r-x)n}_1\over { [t-x]_{n+1,\mathcal{R}(p,q)}}}, 
			\end{equation*} 
			where $t> \kappa + x.$
		\end{lemma}
		{\it Proof:} Setting $v={\epsilon^{t-r-x}_1\over [t]_{\mathcal{R}(p,q)}}$ in (\ref{s24}) and (\ref{s25}), we obtain
		\begin{equation*}
		\sum_{n=\kappa}^{\infty}S_{\mathcal{R}(p,q)}(n;\kappa,x)v^n=\sum_{n=\kappa}^{\infty}S_{\mathcal{R}(p,q)}(n;\kappa,x)\epsilon^{(t-r-x)n}_1[t]^{-n}_{\mathcal{R}(p,q)}
		\end{equation*}
		and 
		\begin{equation*}
		\prod_{j=0}^{\kappa}\Big(1 - \epsilon^{x-j}_1[x+j]_{\mathcal{R}(p,q)}\,
		\epsilon^{t-r-x}_1[t]^{-1}_{\mathcal{R}(p,q)}\Big)^{-1}
		={\epsilon^{(t-r-x)\kappa}_1\,[t]_{\mathcal{R}(p,q)}\over 	\epsilon^{(t-r-x)\kappa}_1[t]^{-\kappa}_{\mathcal{R}(p,q)}}{\epsilon^{-{\kappa+1\choose 2}-x(\kappa+1)}_2\over {[t-x]_{\kappa+1,\mathcal{R}(p,q)}}}. 
		\end{equation*}
		Thus, \begin{equation*}
		\sum_{n=\kappa}^{\infty}S_{\mathcal{R}(p,q)}(n;\kappa,x)\epsilon^{(t-r-x)n}_1[t]^{-n}_{\mathcal{R}(p,q)}={\epsilon^{(t-r-x)\kappa}_1\,[t]_{\mathcal{R}(p,q)}\over \epsilon^{{\kappa+1\choose 2}+x(\kappa+1)}_2\,{[t-x]_{\kappa+1,\mathcal{R}(p,q)}}
		},
		\end{equation*}
		and after rearranging, we find (\ref{s31}). Moreover, let us fix $\kappa$ in (\ref{s31}). Replacing $n$ by $m$ and $\kappa$ by $n,$ we get
		\begin{equation*}
		\epsilon^{{n+1\choose 2}+x(n+1)}_2\sum_{m=n}^{\infty}S_{\mathcal{R}(p,q)}(m;n,x){\epsilon^{(t-r-x)m}_1\over [t]^{m+1}_{\mathcal{R}(p,q)}}={\epsilon^{(t-r-x)n}_1\over {[t-x]_{n+1,\mathcal{R}(p,q)}}}.
		\end{equation*}
		Multiplying the result by $\epsilon^{-{n+1\choose 2}-x(n+1)}_2\,s_{\mathcal{R}(p,q)}(n;\kappa,x),$ and 
		summing for $n=\kappa,\kappa+1,\cdots,$ we find
		\begin{eqnarray*}
			\sum_{n=\kappa}^{\infty}\sum_{m=n}^{\infty}s_{\mathcal{R}(p,q)}(n;\kappa,x)S_{\mathcal{R}(p,q)}(m;n,x){\epsilon^{(t-r-x)m}_1\over [t]^{m+1}_{\mathcal{R}(p,q)}}&=&\sum_{n=\kappa}^{\infty}\epsilon^{{-{n+1\choose 2}}-x(n+1)}_2\nonumber\\
			\quad\quad\quad&\times&s_{\mathcal{R}(p,q)}(n;\kappa,x){\epsilon^{(t-r-x)n}_1\over {[t-x]_{n+1,\mathcal{R}(p,q)}}} 
		\end{eqnarray*}
		By the orthogonality relation (\ref{s14}), we have
		\begin{eqnarray*}
			\sum_{m=\kappa}^{\infty}\delta_{m,\kappa}{\epsilon^{(t-r-x)m}_1\over [t]^{m+1}_{\mathcal{R}(p,q)}}&=&\sum_{n=\kappa}^{\infty}\epsilon^{{-{n+1\choose 2}}-x(n+1)}_2\,s_{\mathcal{R}(p,q)}(n;\kappa,x){\epsilon^{(t-r-x)n}_1\over {[t-x]_{n+1,\mathcal{R}(p,q)}}} 
		\end{eqnarray*}
		and the result follows. $\cqfd$
		\begin{theorem}
			The noncentral $\mathcal{R}(p,q)-$ deformed Stirling numbers of the first and second kinds are given, respectively, by
			\begin{equation}\label{s48}
			s_{\mathcal{R}(p,q)}(n,\kappa,r)= {\epsilon^{-{n\choose 2}+nx}_1\over (\epsilon_1-\epsilon_2)^{n-\kappa}}\sum_{j=\kappa}^{n}(-1)^{j-\kappa}{\epsilon^{{j\choose 2}-r(n-j)-\kappa\,x}_1\over \epsilon^{-{n-j\choose 2}-r(n-j)}_2}\,\bigg[\begin{array}{c} n  \\ j\end{array} \bigg]_{\mathcal{R}(p,q)}\,{j\choose \kappa}
			\end{equation}
			and
			\begin{equation}\label{s49}
			S_{\mathcal{R}(p,q)}(n,\kappa,r)= {\epsilon^{nx}_1 }\sum_{j=\kappa}^{n}(-1)^{j-\kappa}\,{\epsilon^{(n-j)r}_1\,\epsilon^{{\kappa\choose 2}- \kappa\,x}_1\over(\epsilon_1-\epsilon_2)^{n-\kappa}\,\epsilon^{-r(j-\kappa)}_2 }\,{n\choose j}\bigg[\begin{array}{c} j  \\ \kappa \end{array} \bigg]_{\mathcal{R}(p,q)},
			\end{equation}
			where $n\in\mathbb{N}\backslash\{0\},$ $\kappa\in\{1,2,\cdots, n\},$ and $x\in\mathbb{N}.$
		\end{theorem}
		{\it Proof:} From  the relation (\ref{016}) and replacing $p=p^{-1},$ $q=q^{-1},$ $\epsilon_1=\epsilon^{-1}_1,$ $\epsilon_2=\epsilon^{-1}_2,$ $x=\epsilon^{x-r}_1,$ $y=-\epsilon^{x-r}_2,$ $\kappa=j,$ $n=\kappa$ in the $\mathcal{R}(p,q)-$ binomial formula (\ref{bn1}), we get
		\begin{equation}\label{s50}
		{{[x-r]_{n,\mathcal{R}(p,q)}}\over (\epsilon_1-\epsilon_2)^{-n}}=\sum_{j=0}^{n}(-1)^j{\epsilon^{-{n-j\choose 2}+(n-j)(x-r)}_1\over \epsilon^{{j\choose 2}-j(x-r)}_2}(\epsilon_1\,\epsilon_2)^{-j(n-j)}
		\bigg[\begin{array}{c} n  \\ j \end{array} \bigg]_{\mathcal{R}(p,q)}.
		\end{equation}
		Multipliying (\ref{s50}) by $\epsilon^{{n\choose 2}+rn}_2,$ and using
		\begin{eqnarray*}
			\epsilon^{-jx}_1\,\epsilon^{jx}_2 
			&=& \sum_{\kappa=0}^{j}(-1)^\kappa\,{j\choose \kappa}\,\epsilon^{-\kappa\,x}_1\,(\epsilon_1-\epsilon_2)^\kappa\,{[x]^\kappa_{\mathcal{R}(p,q)}},
		\end{eqnarray*}
		we obtain
		\begin{eqnarray*}
			&&{{[x-r]_{n,\mathcal{R}(p,q)}}\over \epsilon^{-{n\choose 2}-rn}_2}=\sum_{\kappa=0}^{n}\bigg({ \epsilon^{-{n\choose 2}+nx}_1\over (\epsilon_1-\epsilon_2)^{n-\kappa}}\sum_{j=\kappa}^{n}(-1)^{j-\kappa} \epsilon^{{j\choose 2}-r(n-j)-\kappa\,x}_1\cr&&
			\qquad\qquad\qquad\qquad\qquad\times \epsilon^{{n-j\choose 2}+r(n-j)}_2
			\bigg[\begin{array}{c} n  \\ j \end{array} \bigg]_{\mathcal{R}(p,q)}\,{j\choose \kappa}\bigg){[x]^\kappa_{\mathcal{R}(p,q)}}.
		\end{eqnarray*}
		Comparing the above relation with (\ref{s4}), we obtain (\ref{s48}). Furthermore,
		\begin{eqnarray}\label{s55}
		[x+r]^n_{\mathcal{R}(p,q)}
		&=& { \epsilon^{nx}_1\over (\epsilon_1-\epsilon_2)^n}\sum_{j=0}^{n}(-1)^j\,{n\choose j}\,\epsilon^{(n-j)r}_1\,\epsilon^{jr}_2\,\epsilon^{-jx}_1\,\epsilon^{jx}_2.
		\end{eqnarray}
		Using (\ref{bna11}), we get 
		\begin{eqnarray*}
			[x+r]^n_{\mathcal{R}(p,q)}&=&{ \epsilon^{nx}_1\over (\epsilon_1-\epsilon_2)^n}\sum_{j=0}^{n}(-1)^j\,{n\choose j}\,\epsilon^{(n-j)r}_1\,\epsilon^{jr}_2\nonumber\\
			&\times& \sum_{\kappa=0}^{j}(-1)^{\kappa}\,\epsilon_1^{{\kappa \choose 2}-\kappa \,x}\,\epsilon_2^{{\kappa \choose 2}} \bigg[\begin{array}{c} j   \\ \kappa \end{array} \bigg]_{\mathcal{R}(p,q)}(\epsilon_1-\epsilon_2)^\kappa\,	[x]_{\kappa,\mathcal{R}(p,q)},
		\end{eqnarray*}
		and from (\ref{s6}), we obtain (\ref{s49}).
		$\cqfd$
		\begin{remark}
			Taking $\mathcal{R}(p,q)=q,$ we obtain the  noncentral $q-$  Stirling numbers of the first and second kinds as:
			\begin{equation*}\label{qs48}
			s_{q}(n,\kappa,r)= {1\over (1-q)^{n-\kappa}}\sum_{j=\kappa}^{n}(-1)^{j-\kappa}{ q^{{n-j\choose 2}+r(n-j)}}\,\bigg[\begin{array}{c} n  \\ j\end{array} \bigg]_{q}\,{j\choose \kappa}
			\end{equation*}
			and
			\begin{equation*}\label{qs49}
			S_{q}(n,\kappa,r)= {1\over(1-q)^{n-\kappa}}\sum_{j=\kappa}^{n}(-1)^{j-\kappa}\,q^{r(j-\kappa)} \,{n\choose j}\bigg[\begin{array}{c} j  \\ \kappa \end{array} \bigg]_{q},
			\end{equation*}
			where $n\in\mathbb{N}\backslash\{0\},$ $\kappa\in\{1,2,\cdots, n\},$ and $x\in\mathbb{N}.$
		\end{remark}
		\begin{corollary}
			Let $\kappa$ and $j$ be positive integers. Then, the following relations hold:
			\begin{equation}\label{s33}
			{\kappa \choose j} =\sum_{m=j}^{\kappa}(-1)^{m-j}(\epsilon_1-\epsilon_2)^{m-j} \epsilon^{{m\choose 2}-x(m-j)}_1 s_{\mathcal{R}(p,q)}(m,j)\,\bigg[\begin{array}{c} \kappa  \\ m \end{array} \bigg]_{\mathcal{R}(p,q)}
			\end{equation} 
			and
			\begin{equation}\label{s33b}
			\bigg[\begin{array}{c} \kappa  \\ j \end{array} \bigg]_{\mathcal{R}(p,q)}=\sum_{m=j}^{\kappa}(-1)^{m-j}(\epsilon_1-\epsilon_2)^{m-j} \epsilon^{-{j\choose 2}-x(m-j)}_1 S_{\mathcal{R}(p,q)}(m,j)\,{\kappa \choose m},
			\end{equation}
			where $x$ is an integer and $0<q<p\leq 1.$
		\end{corollary}
		{\it Proof:} Replacing $j$ by $i,$ $n$ by 
		$m,$ $\kappa$ by $j$ and $r=0$ in
		(\ref{s48}), we get
		\begin{equation*}
		s_{\mathcal{R}(p,q)}(m,j)= {\epsilon^{-{m\choose 2}+m\,x}_1\over (\epsilon_1-\epsilon_2)^{m-j}}\sum_{i=j}^{m}(-1)^{i-j}\epsilon^{{i\choose 2}-j\,x}_1\epsilon^{{m-i\choose 2}}_2\,\bigg[\begin{array}{c} m  \\ i\end{array} \bigg]_{\mathcal{R}(p,q)}\,{i\choose j}.
		\end{equation*}
		Multiplying this result by $$(-1)^{m-j}{ (\epsilon_1-\epsilon_2)^{m-j}\over \epsilon^{-{m\choose 2}+(m-j)x}_1 }\bigg[\begin{array}{c} \kappa  \\ m \end{array} \bigg]_{\mathcal{R}(p,q)}$$ and summing for all $m=j,j+1,\cdots,\kappa,$ we obtain 
		\begin{equation*}
		\sum_{m=j}^{\kappa}(-1)^{m-j}{ (\epsilon_1-\epsilon_2)^{m-j}\over \epsilon^{-{m\choose 2}+(m-j)x}_1 }s_{\mathcal{R}(p,q)}(m,j)\bigg[\begin{array}{c} \kappa  \\ m \end{array} \bigg]_{\mathcal{R}(p,q)}
		= \sum_{i=j}^{\kappa}{i \choose j}\,\delta_{\kappa,i}={\kappa \choose j}.
		\end{equation*}
		Similarly, in (\ref{s49}), we replace $n$ by $m,$ $j$ by $i,$ $r=0,$   multiply the resulting expression by
		$$(-1)^{m-j}\,{ (\epsilon_1-\epsilon_2)^{m-j}\over \epsilon^{{j\choose 2}+(m-j)x}_1 }\,{ \kappa  \choose m },$$ and sum it for all $m=j,j+1,\cdots, \kappa,$ to get the result.
		$\cqfd$
		
		Note that we obtain the $q-$ deformed formulae (\ref{qs33}) and (\ref{qs33b}) by taking $\mathcal{R}(x,y)=y.$
		\begin{lemma}
			Let $u$ be a natural number. Then, the following relations hold:
			\begin{equation}\label{sf1}
			s_{\mathcal{R}(p,q)}(u,1)=(-1)^{u-1}\,\epsilon^{xu-{(u-1)(u+2)\over 2}}_1\mathcal{R}!(p^{u-1},q^{u-1})
			\end{equation}
			and
			\begin{equation*}\label{sf2}
			s_{\mathcal{R}(p,q)}(u,2)=(-1)^{u-2}\,\epsilon^{x(u-2)-{(u-2)(u+3)\over 2}}_1\,\mathcal{R}!(p^{u-1},q^{u-1})\zeta_{u-1,p,q}
			\end{equation*}
			where $\zeta_{u,p,q}=\displaystyle\sum_{j=1}^{u}\epsilon^{(j-1)(x+1)}_1\Bigg(\mathcal{R}(p^j,q^j)\Bigg)^{-1}$ and $x\in\mathbb{N}.$
		\end{lemma}
		{\it Proof:} From the triangular recursion relation of the $\mathcal{R}(p,q)-$ deformed Stirling numbers of the first kind (\ref{s16}), we set $\kappa=1,$ $j=0,$ and obtain the first-order recursion relation
		\begin{equation*}
		s_{\mathcal{R}(p,q)}(n,1)=-\epsilon^{x-n}_1\,[n-1]_{\mathcal{R}(p,q)}s_{\mathcal{R}(p,q)}(n-1,1)
		\end{equation*}
		with $n\in\mathbb{N}\backslash\{0,1\}$ and  $s_{\mathcal{R}(p,q)}(1,1)=s_{\mathcal{R}(p,q)}(0,0):=1.$ By iteration, we get (\ref{sf1}). Setting  $\kappa=2$ also leads to the recursion relation
		\begin{equation}
		s_{\mathcal{R}(p,q)}(n,2)+\epsilon^{x-n}_1\,[n-1]_{\mathcal{R}(p,q)}s_{\mathcal{R}(p,q)}(n-1,2)=R(n,p,q),
		\end{equation}
		where
		$$R(n,p,q)=(-1)^{n-2}[n-2]_{\mathcal{R}(p,q)}!\,\epsilon^{x(n-1)-{(n-2)(n+1)\over 2}}_1,$$ 
		which is solved to give the required expression. $\cqfd$
		
		Taking $\mathcal{R}(p,q)=q,$ we recover
		\begin{equation}\label{qsf1}
		s_{q}(u,1)=(-1)^{u-1}\,[u-1]_q!
		\end{equation}
		and
		\begin{equation}\label{qsf2}
		s_{q}(u,2)=(-1)^{u-2}[u-1]_q!\zeta_{u-1,q},
		\end{equation}
		where $\zeta_{u,p,q}=\displaystyle\sum_{j=1}^{u}[j]^{-1}_q!.$
		\begin{remark} The particular case of the generalized $q-$Quesne quantum algebra is here worthy of attention as matter of illustration.
			%
			Indeed, the signless (or absolute) noncentral \textbf{generalized $q-$ Quesne}  Stirling number of the
			first kind is given by 
			\begin{equation*}\label{s7g}
			|s^Q_{p^{-1},q^{-1}}(n, \kappa,j)| = \Big({q\over p}\,[-1]^Q_{p,q}\Big)^{n-\kappa}\,\,s^Q_{p^{-1},q^{-1}}(n, \kappa,j)
			\end{equation*}
			or 
			\begin{equation*}\label{s8g}
			|s^Q_{p,q}(n, \kappa,j)| = \Big({q\over p}\,[-1]^Q_{p,q}\Big)^{-n+\kappa}\,\,s^Q_{p,q}(n, \kappa,j),
			\end{equation*}
			where $n$ is a positive interger.
			The  related deformed Stirling numbers of the first and second kinds verify the orthogonality relations
			\begin{equation*}\label{s13g}
			\sum_{m=\kappa}^{n}s^Q_{p,q}(n,m,j)\,S^Q_{p,q}(m,\kappa,j)= \delta_{n,\kappa}
			\end{equation*}
			and 
			\begin{equation*}\label{s14g}
			\sum_{m=\kappa}^{n}S^Q_{p,q}(n,m,j)\,s^Q_{p,q}(m,\kappa,j)= \delta_{n,\kappa}.
			\end{equation*}
			They  obey, respectively, 
			the following recursion relations:
			\begin{equation*}\label{s16g}
			s^Q_{p,q}(n+1,\kappa;j)= s^Q_{p,q}(n,\kappa-1;j)- q\,p^{x-n-j-1}[n+j]^Q_{p,q}s^Q_{p,q}(n,\kappa,j)
			\end{equation*}
			with the initial conditions
			$s^Q_{p,q}(0,0,j)=1,$  $s^Q_{p,q}(n,0,j)={q^{-{n \choose 2}-n(j-1)}\over p^n}\,{[-j]^Q_{n,p,q}},$
			$n\geq 0.$ and $s^Q_{p,q}(0,\kappa,j)=0,$ $j\geq 0;$ and
			\begin{equation*}\label{s17g}
			S^Q_{p,q}(n+1,\kappa;j)= p^{j}\,S^Q_{p,q}(n,\kappa-1;j)+ {q\over p}\,[\kappa+j]^Q_{p,q}\,S^Q_{p,q}(n,\kappa;j)
			\end{equation*}
			with the initial conditions $s^Q_{p,q}(0,0,j)=1,$  $s^Q_{p,q}(n,0,j)=\big({q\over p}[j]^Q_{p,q}\big)^n,$
			$n\geq 0,$  and  $s^Q_{p,q}(0,\kappa,j)=0,$ $\kappa\geq 0.$ 
			For fix $\kappa$, they are generated by the function given by
			\begin{equation*}\label{s24g}
			\Psi_{\kappa}(v;p,q,r) = \sum_{n=\kappa}^{\infty}S^Q_{p,q}(n;\kappa,r)\,v^n\mbox{,}\quad \kappa\in\mathbb{N}
			\end{equation*}
			developed in the product form as:
			\begin{equation*}\label{s25g}
			\Psi_{\kappa}(v;p,q,r) = v^{\kappa}\,\prod_{j=0}^{\kappa}\Big(1 - q\,p^{x-j-1}\,[r+j]^Q_{p,q}\,v\Big)^{-1}
			\end{equation*}
			for $|v|< {p^{-x+j+1} \over q\,[r+j]^Q_{p,q}}.$
			Their reciprocal
			factorial ${[t-x]^Q_{-\kappa-1,p,q}}$ is expanded into reciprocal generalized $q-$ Quesne powers $\big([t]^Q_{p,q}\big)^{-n-1}$ as follows:
			\begin{equation*}\label{s31g}
			{p^{(t-r-x+1)\kappa+1}\over q\,{[t-x]^Q_{\kappa+1,p,q}}} = q^{-{\kappa +1\choose 2}-x(\kappa+1)}\sum_{n=\kappa}^{\infty}S^Q_{p,q}(n,\kappa;x){p^{(t-r-x+1)n+1}\over \big(q\,[t]^Q_{p,q}\big)^{n+1}},
			\end{equation*} 
			while  their reciprocal 
			powers $\big([t]_{p,q}^Q\big)^{-\kappa-1}$  are spanned in   the reciprocal noncentral generalized $q-$ Quesne  factorial ${[t-x]_{-n-1,p,q}^Q}$ as:
			\begin{equation*}\label{s32g}
			{p^{(t-r-x+1)\kappa+1}\over \big(q\,[t]_{p,q}^Q\big)^{\kappa+1}} =\sum_{n=\kappa}^{\infty}\,q^{-{n +1\choose 2}+x(n+1)}s^Q_{p,q}(n,\kappa;x){p^{(t-r-x+1)n+1}\over {q\,[t-x]^Q_{n+1,p,q}}}, 
			\end{equation*} 
			where $t> \kappa + x.$
			Moreover,
			\begin{equation*}\label{s48g}
			s^Q_{p,q}(n,\kappa,r)= {p^{-{n\choose 2}}\over p^{-nx} }\sum_{j=\kappa}^{n}(-1)^{j-\kappa}\,{p^{{j\choose 2}-r(n-j)-x\kappa}{j\choose \kappa}\over (p-q^{-1})^{n-\kappa}\,q^{{n-j\choose 2}+r(n-j)}}\,\bigg[\begin{array}{c} n  \\ j\end{array} \bigg]^Q_{p,q}
			\end{equation*}
			and
			\begin{equation*}\label{s49g}
			S^Q_{p,q}(n,\kappa,r)= p^{n\,x}\sum_{j=\kappa}^{n}(-1)^{j-\kappa}\,{p^{(n-j)r}\,p^{{\kappa\choose 2}-\kappa\,x}\over (p-q^{-1})^{n-\kappa}\,q^{r(j-\kappa)}}{n\choose j}\,\bigg[\begin{array}{c} j  \\ \kappa \end{array} \bigg]^Q_{p,q},
			\end{equation*}
			where $n\in\mathbb{N}\backslash\{0\},$ $\kappa\in\{1,2,\cdots, n\},$ and $x\in\mathbb{N}.$
			In particular,
			\begin{equation*}
			S^Q_{p,q}(n,\kappa,r)={(q^{-1}\,p)^{-1+\kappa}p^{nx}\over [\kappa]^Q_{p,q}!}\sum_{i=0}^{\kappa}(-1)^{\kappa-i}\,{p^{{i\choose 2}+j(r+i)}\over q^{{i+1\choose 2}-\kappa(r+i)}}\bigg[\begin{array}{c} \kappa  \\ i \end{array} \bigg]^Q_{p,q}\,[r+i]^Q_{p,q}.
			\end{equation*}
			For  positive integers $\kappa$ and $j,$
			\begin{equation*}
			{\kappa \choose j} =\sum_{m=j}^{\kappa}(-1)^{m-j}\,{(p-q^{-1})^{m-j}\over p^{-{m\choose 2}+x(m-j)}}\, s^Q_{p,q}(m,j)\,\bigg[\begin{array}{c} \kappa  \\ m \end{array} \bigg]^Q_{p,q}
			\end{equation*} 
			and
			\begin{equation*}
			\bigg[\begin{array}{c} \kappa  \\ j \end{array} \bigg]^Q_{p,q}=\sum_{m=j}^{\kappa}(-1)^{m-j}\,{(p-q^{-1})^{m-j} \over p^{{j\choose 2}+x(m-j)}}\, S^Q_{p,q}(m,j)\,{\kappa \choose m},
			\end{equation*}
			where $x$ is an integer and $0<q<p\leq 1.$
			Furthermore, given   a natural number $u$, 
			\begin{equation*}\label{sfg1}
			s_{\mathcal{R}(p,q)}(u,1)=(-1)^{u-1}{p^{2xu-u^2+u }\over 2q^{u-1}}\,[u-1]^Q_{p,q}!
			\end{equation*}
			and
			\begin{equation*}\label{sfg2}
			s_{p,q^{-1}}(u,2)=(-1)^{u-2}{p^{2x(u-2)-u^2+u+4}\over 2q^{u-1}}\,[u-1]^Q_{p,q}!\zeta_{u-1,p,q},
			\end{equation*}
			where $\zeta_{u,p,q}=\displaystyle\sum_{j=1}^{u}{p^{(j-1)(x+1)+1}\over q\,[j]^Q_{p,q}}.$	
		\end{remark}
		\subsection{$\mathcal{R}(p,q)-$ deformed  Bell numbers}
		Let us consider a simple finite graph $\bf G$ with $n$ vertices, and $\kappa$   independent blocks. $V(\bf G)$ is a set of vertices of $\bf G.$ The vertices are partitioned into $\kappa$ independent blocks by $\Pi=\{V_1,V_2,\cdots,V_{\kappa}\}.$ We define the $\mathcal{R}(p,q)-$deformed weight as follows:
		\begin{equation*}
		W_{\mathcal{R}(p,q)}(\pi) :=\Big({\epsilon_2\over \epsilon_1}\Big)^{\displaystyle\sum_{j=1}^{\kappa}(j-1)|V_j|},
		\end{equation*}
		where $|V_j|$ is the cardinality of the set $V.$ The $\mathcal{R}(p,q)-$ deformed Stirling number of the second kind for the graph is expressed by
		\begin{equation*}
		S_{\mathcal{R}(p,q)}\big({\bf G}, \kappa\big):= \sum_{\Lambda} W_{\mathcal{R}(p,q)}(\pi),
		\end{equation*}
		where $S_{\mathcal{R}(p,q)}\big(\bf G, 0\big):=0,$ and $\Lambda$ denotes the independent partitions $\pi$ of $V(\bf G).$ Analogously, the $\mathcal{R}(p,q)-$ deformed Bell number for the graph is defined by
		\begin{equation*}
		B_{\mathcal{R}(p,q)}\big( G\big):= \sum_{\kappa=0}^{|V(\bf G)|}S_{\mathcal{R}(p,q)}\big(\bf G, \kappa\big).
		\end{equation*}
		\begin{theorem}
			Let $\bar{T}_n$ be the dual path graph of $\bf G.$ Then, the $\mathcal{R}(p,q)-$ deformed Stirling numbers of the second kind and the $\mathcal{R}(p,q)-$ deformed Bell numbers for the graph $\bf G$ are given, respectively, by 
			\begin{equation*}
			S_{\mathcal{R}(p,q)}\big(\bar{T}_n, \kappa\big) = {\epsilon^{{n\choose 2}-\kappa(n-\kappa)}_2\over \epsilon^{{n\choose 2}-\kappa(n-\kappa)+\kappa-1}_1}\,\bigg[\begin{array}{c} \kappa  \\ n-\kappa \end{array} \bigg]_{\mathcal{R}(p,q)}
			\end{equation*}
			and
			\begin{equation*}
			B_{\mathcal{R}(p,q)}\big(\bar{T}_n\big)=\sum_{\kappa=0}^{n} {\epsilon^{{n\choose 2}-\kappa(n-\kappa)}_2\over \epsilon^{{n\choose 2}-\kappa(n-\kappa)+\kappa-1}_1}\,\bigg[\begin{array}{c} \kappa  \\ n-\kappa \end{array} \bigg]_{\mathcal{R}(p,q)},
			\end{equation*}
			where $n\in\mathbb{N}.$
		\end{theorem}
		{\it Proof:} Let us consider the expression
		\begin{equation}\label{bn}
		S_{\mathcal{R}(p,q)}(\bar{T}_n, \kappa) = {\epsilon^{{\kappa\choose 2}+ {n-\kappa\choose 2}}_2\over \epsilon^{{\kappa\choose 2}+ {n-\kappa\choose 2}+\kappa-1}_1}\,\bigg[\begin{array}{c} \kappa  \\ n-\kappa \end{array} \bigg]_{\mathcal{R}(p,q)}.
		\end{equation}
		For $n=0,$ the relation (\ref{bn}) is true. We assume that (\ref{bn}) is true for all $n$ and prove it for $l=n+1.$ Consider
		\begin{equation*}
		S_{\mathcal{R}(p,q)}(\bar{T}_{n+1}, \kappa) =\hat{S}_{\mathcal{R}(p,q)}(\bar{T}_{n+1}, \kappa) + \tilde{S}_{\mathcal{R}(p,q)}(\bar{T}_{n+1}, \kappa),
		\end{equation*}
		where
		\begin{eqnarray*}
			\hat{S}_{\mathcal{R}(p,q)}(\bar{T}_{n+1}, \kappa)&=&{\epsilon^{n-2\kappa+1}_1\over \epsilon^{-\kappa+1}_2}\,S_{\mathcal{R}(p,q)}(\bar{T}_{n}, \kappa-1)\nonumber\\
			&=&{\epsilon^{{\kappa\choose 2}+ {n-\kappa+1\choose 2}}_2\over \epsilon^{{\kappa\choose 2}+ {n-\kappa+1\choose 2}+\kappa-1}_1}\,{\epsilon^{n-\kappa+1}_1}\,\bigg[\begin{array}{c} \kappa-1  \\ n-\kappa+1 \end{array} \bigg]_{\mathcal{R}(p,q)}
		\end{eqnarray*}
		and
		\begin{eqnarray*}
			\tilde{S}_{\mathcal{R}(p,q)}(\bar{T}_{n+1}, \kappa)&=& {\epsilon^{2(\kappa-1)}_2\over \epsilon^n_1}\,S_{\mathcal{R}(p,q)}(\bar{T}_{n-1}, \kappa-1)\nonumber\\
			&=&{\epsilon^{{\kappa\choose 2}+ {n-\kappa+1\choose 2}}_2\over \epsilon^{{\kappa\choose 2}+ {n-\kappa+1\choose 2}+\kappa-1}_1}\,{\epsilon^{2\kappa-1-n}_2}\,\bigg[\begin{array}{c} \kappa-1  \\ n-\kappa \end{array} \bigg]_{\mathcal{R}(p,q)}.
		\end{eqnarray*}
		Using Eq.(\ref{bc1}) and after computation, we obtain 
		\begin{equation*}
		S_{\mathcal{R}(p,q)}(\bar{T}_{n+1}, \kappa)={\epsilon^{{\kappa\choose 2}+ {n-\kappa+1\choose 2}}_2\over \epsilon^{{\kappa\choose 2}+ {n-\kappa+1\choose 2}+\kappa-1}_1}\bigg[\begin{array}{c} \kappa  \\ n-\kappa+1 \end{array} \bigg]_{\mathcal{R}(p,q)}.
		\end{equation*}
		Thus, the proof  is achieved. $\cqfd$
		\subsection{Application
		}
		We consider the dual path graph $\bar{T}_5.$ It has $4$ independent partitions into $4$ blocks given as follows:
		\begin{equation*}
		\Gamma_1=\Big(\{1,2\},\{3\},\{4\},\{5\}\Big),\quad \Gamma_2=\Big(\{1\},\{2,3\},\{4\},\{5\}\Big)
		\end{equation*}
		\begin{equation*}
		\Gamma_3=\Big(\{1\},\{2\},\{3,4\},\{5\}\Big),\quad \Gamma_4=\Big(\{1\},\{2\},\{3\},\{4,5\}\Big).
		\end{equation*}
		The $\mathcal{R}(p,q)-$deformed weight is given by
		\begin{equation*}
		W_{\mathcal{R}(p,q)}(\Gamma) =\Bigg({\epsilon_2\over \epsilon_1}\Bigg)^{\displaystyle\sum_{j=1}^{4}(j-1)|V_j|}.
		\end{equation*}
		Hence,
		\begin{eqnarray*}
			W_{\mathcal{R}(p,q)}(\Gamma_1) &=&\Bigg({\epsilon_2\over \epsilon_1}\Bigg)^{\displaystyle\sum_{j=1}^{4}(j-1)|V_j|}
			=\Bigg ({\epsilon_2\over \epsilon_1}\Bigg)^6.
		\end{eqnarray*}
		Similarly, we get
		$
		W_{\mathcal{R}(p,q)}(\Gamma_2)=\Big({\epsilon_2\over \epsilon_1}\Big)^7,$ $W_{\mathcal{R}(p,q)}(\Gamma_3)=\Big({\epsilon_2\over \epsilon_1}\Big)^8$ and $W_{\mathcal{R}(p,q)}(\Gamma_4)=\Big({\epsilon_2\over \epsilon_1}\Big)^9.$ Finally,
		\begin{eqnarray}
		S_{\mathcal{R}(p,q)}\bigg(\bar{T}_5, 4\bigg)&=&\Big({\epsilon_2\over \epsilon_1}\Big)^6\,{1-({\epsilon_2\over \epsilon_1})^4\over 1-({\epsilon_2\over \epsilon_1})}\nonumber\\
		&=& \Big({\epsilon_2\over \epsilon_1}\Big)^6\,\,\epsilon^{-3}_1\,\mathcal{R}(p^4,q^4)\nonumber\\
		&=&\Big({\epsilon_2\over \epsilon_1}\Big)^6\,\,\epsilon^{-3}_1\,\bigg[\begin{array}{c} 4  \\ 1 \end{array} \bigg]_{\mathcal{R}(p,q)}.
		\end{eqnarray}
		\begin{remark}
			\begin{enumerate}
				\item[(1)] Note that  the $q-$ Stirling number of the second kind and $q-$ Bell number of the graph can easily  be derived by taking  $\mathcal{R}(p,q)=q$ as follows:
				\begin{equation*}
				S_{q}\big(\bar{T}_n, \kappa\big) = {q^{{n\choose 2}-\kappa(n-\kappa)}}\,\bigg[\begin{array}{c} \kappa  \\ n-\kappa \end{array} \bigg]_{q}
				\end{equation*}
				and
				\begin{equation*}
				B_{q}\big(\bar{T}_n\big)=\sum_{\kappa=0}^{n} {q^{{n\choose 2}-\kappa(n-\kappa)}}\,\bigg[\begin{array}{c} \kappa  \\ n-\kappa \end{array} \bigg]_{q},
				\end{equation*}
				where $n\in\mathbb{N}.$
				\item[(2)] The \textbf{generalized $q-$Quesne} Stirling number of the second kind and the Bell number for the dual path graph are, respectively, given by
				\begin{equation*}
				S^Q_{p,q}\big(\bar{T}_n, \kappa\big) = {q^{-{n\choose 2}+\kappa(n-\kappa)}\over p^{{n\choose 2}-\kappa(n-\kappa)+\kappa-1}}\,\bigg[\begin{array}{c} \kappa  \\ n-\kappa \end{array} \bigg]^Q_{p,q}
				\end{equation*}
				and
				\begin{equation*}
				B^Q_{p,q}\big(\bar{T}_n\big)=\sum_{\kappa=0}^{n} {q^{-{n\choose 2}+\kappa(n-\kappa)}\over p^{{n\choose 2}-\kappa(n-\kappa)+\kappa-1}}\,\bigg[\begin{array}{c} \kappa  \\ n-\kappa \end{array} \bigg]^Q_{p,q},
				\end{equation*}
				where $n\in\mathbb{N}.$
			\end{enumerate}
		\end{remark}
		\subsection{$\mathcal{R}(p,q)-$ deformed factorial and binomial moments}
		For a study on  $q-$ factorial and $q-$ binomial moments, see \cite{CA1}. We deal here with the $\mathcal{R}(p,q)-$generalization. For that, we consider 
		a nonnegative integer-valued discrete random variable $X,$ and $g(x)=P(X=x),$ $x\in\mathbb{N}\backslash\{0\},$ the probability distribution of $X.$  
		\begin{equation*}
		{\bf E}\Big({[X]_{r,\mathcal{R}(p,q)}}\Big) = \sum_{x=r}^{\infty}{[x]_{r,\mathcal{R}(p,q)}}\,g(x)\mbox{,}\quad r\in\mathbb{N}\backslash\{0\},
		\end{equation*}
		\begin{equation}\label{fb2}
		{\bf E}\bigg(\bigg[\begin{array}{c} X  \\ r\end{array} \bigg]_{\mathcal{R}(p,q)}\bigg)  = \sum_{x=r}^{\infty}\bigg[\begin{array}{c} x  \\ r\end{array} \bigg]_{\mathcal{R}(p,q)}\,g(x)\mbox{,}\quad r\in\mathbb{N}\backslash\{0\}
		\end{equation}
		referred to the $r^{th}-$order $\mathcal{R}(p,q)-$ factorial and  $r^{th}-$order $\mathcal{R}(p,q)-$ binomial moments, respectively, of the random variable $X.$
		In the particular case of $r=1,$ we define the $\mathcal{R}(p,q)-$ mean value, also called the $\mathcal{R}(p,q)-$ expectation value, of $X$  by 
		\begin{equation*}
		\mu_{\mathcal{R}(p,q)}:={\bf E}\Big(\big[X\big]_{\mathcal{R}(p,q)}\Big)=\sum_{x=1}^{\infty}\big[x\big]_{\mathcal{R}(p,q)}\,g(x).
		\end{equation*}
		The $\mathcal{R}(p,q)-$ variance of $X$ is then obtained as
		\begin{equation*}
		\sigma^2_{\mathcal{R}(p,q)}:={\bf V}\Big(\big[X\big]_{\mathcal{R}(p,q)}\Big)= {\bf E}\big(\big[X\big]_{\mathcal{R}(p,q)}\big)^2 - \bigg[{\bf E}\Big(\big[X\big]_{\mathcal{R}(p,q)}\Big)\bigg]^2.
		\end{equation*}
		Since   $[X-1]_{\mathcal{R}(p,q)}=\epsilon^{-1}_2[X]_{\mathcal{R}(p,q)}-\epsilon^{-1}_2\epsilon^{X-1}_1,$ and ${[X]_{2,\mathcal{R}(p,q)}}= [X]_{\mathcal{R}(p,q)}[X-1]_{\mathcal{R}(p,q)},$  then
		\begin{equation*}
		{\bf V}\Big([X]_{\mathcal{R}(p,q)}\Big)=\epsilon_2\,{\bf E}\big({[X]_{2,\mathcal{R}(p,q)}}\big) +\epsilon^{X-1}_1\,{\bf E}\big([X]_{\mathcal{R}(p,q)}\big)- \big[{\bf E}\big([X]_{\mathcal{R}(p,q)}\big)\big]^2.
		\end{equation*}
		\begin{theorem}
			The binomial moment is given as function of the $\mathcal{R}(p,q)-$ binomial moment as follows:
			\begin{equation}\label{fb6}
			{\bf E}\bigg(\bigg[\begin{array}{c} X  \\ j\end{array} \bigg]_{\mathcal{R}(p,q)}\bigg)= \sum_{m=j}^{\infty} (-1)^{m-j}{(\epsilon_1-\epsilon_2)^{m-j}\over \epsilon^{-{m\choose 2}+\tau(m-j)}_1}s_{\mathcal{R}(p,q)}(m,j){\bf E}\bigg(\bigg[\begin{array}{c} X  \\ m \end{array} \bigg]_{\mathcal{R}(p,q)}\bigg),
			\end{equation}
			while the factorial moment is given in terms of the $\mathcal{R}(p,q)-$ factorial moment by
			\begin{equation}\label{fb7}
			{\bf E}[(X)_j]= j!\sum_{m=j}^{\infty} (-1)^{m-j}{(\epsilon_1-\epsilon_2)^{m-j}\over \epsilon^{-{m\choose 2}+\tau(m-j)}_1}s_{\mathcal{R}(p,q)}(m,j){{\bf E}\big({[ X ]_{m,\mathcal{R}(p,q)}}\big)\over [m]_{\mathcal{R}(p,q)}!},
			\end{equation}
			where $j\in\mathbb{N}\backslash\{0\},$ $\tau\in\mathbb{N},$ and $s_{\mathcal{R}(p,q)}$ is the $\mathcal{R}(p,q)-$ deformed Stirling number of the first kind.
		\end{theorem}
		{\it Proof:} Multiplying (\ref{s33}) by the probability distribution $g(x)$ and summing for all $x\in\mathbb{N},$ we deduce   (\ref{fb2}) from (\ref{fb6}). Moreover, from 
		\begin{equation*}
		{\bf E}\bigg[{X\choose j}\bigg]= { {\bf E}[(X)_j]\over j!},\quad {\bf E}\bigg(\bigg[\begin{array}{c} X  \\ m \end{array} \bigg]_{\mathcal{R}(p,q)}\bigg)={{\bf E}\big({[X]_{m,\mathcal{R}(p,q)}}\big)\over [m]_{\mathcal{R}(p,q)}!}
		\end{equation*} 
		and the relation (\ref{fb6}), we derive (\ref{fb7}). $\cqfd$
		\begin{remark}
			Putting $\mathcal{R}(u,v)=v,$ we obtain  the usual binomial moment as function of the $q-$ binomial moment as:
			\begin{equation*}\label{fba6}
			{\bf E}\bigg(\bigg[\begin{array}{c} X  \\ j\end{array} \bigg]_{q}\bigg)= \sum_{m=j}^{\infty} (-1)^{m-j}{(1-q)^{m-j}}s_{q}(m,j){\bf E}\bigg(\bigg[\begin{array}{c} X  \\ m \end{array} \bigg]_{q}\bigg),
			\end{equation*}
			while the usual factorial moment is given in terms of the $q-$ factorial moment by
			\begin{equation*}\label{fba7}
			{\bf E}[(X)_j]= j!\sum_{m=j}^{\infty} (-1)^{m-j}{(1-q)^{m-j}}s_{q}(m,j){{\bf E}\big({[ X ]_{m,q}}\big)\over [m]_{q}!},
			\end{equation*}
			where $j\in\mathbb{N}\backslash\{0\},$ and $s_{q}$ is the $q-$Stirling number of the first kind.
		\end{remark}
		\begin{theorem}
			The ${\mathcal{R}(p,q)}-$deformed probability distribution $g(x)$ of a discrete random variable $X$ is given by the absolutely convergent series
			\begin{equation*}\label{fb9}
			g(x) = \sum_{m=x}^{\infty}(-1)^{m-x}\epsilon^{{x\choose 2}}_1\,\epsilon^{{m-x\choose 2}}_2\bigg[\begin{array}{c} m \\ x \end{array} \bigg]_{\mathcal{R}(p,q)}{\bf E}\bigg(\bigg[\begin{array}{c} X  \\ m \end{array} \bigg]_{\mathcal{R}(p,q)}\bigg),\quad x\in\mathbb{N}.
			\end{equation*}
		\end{theorem}
		{\it Proof:} Replacing $x$ by $\kappa$ and $r$ by $m$ in expression (\ref{fb2}), multiplying it by $$(-1)^{m-x}\epsilon^{{x\choose 2}}_1\,\epsilon^{{m-x\choose 2}}_2\bigg[\begin{array}{c} m \\ x \end{array} \bigg]_{\mathcal{R}(p,q)},$$ and 
		summing for all $m=x,x+1,\cdots,$ we obtain
		\begin{equation*}
		\sum_{m=x}^{\infty}(-1)^{m-x}\epsilon^{{x\choose 2}}_1\,\epsilon^{{m-x\choose 2}}_2\bigg[\begin{array}{c} m \\ x \end{array} \bigg]_{\mathcal{R}(p,q)}{\bf E}\bigg(\bigg[\begin{array}{c} X  \\ m \end{array} \bigg]_{\mathcal{R}(p,q)}\bigg) 
		=\sum_{\kappa=x}^{\infty}\delta_{\kappa,x}g(\kappa)=g(x).
		\end{equation*}
		$\cqfd$
		
		Note that the  probability distribution (\ref{fba9}) can be retrieved by taking $\mathcal{R}(u,v)=v.$ 
		\begin{remark}The particular case of the \textbf{generalized $q-$Quesne}  factorial and binomial moment,
			and  probability distribution is detailed as follows: 
			\begin{enumerate}
				\item[(1)] Let $X$ be a nonnegative integer-valued discrete random variable and $h(x)=P(X=x),$ $x\in\mathbb{N}\backslash\{0\},$ the probability distribution of $X.$ Assume the convergence  of the series:  
				\begin{equation*}
				{\bf E}\Big({[X]^Q_{r,p,q}}\Big) = \sum_{x=r}^{\infty}{[x]^Q_{r,p,q}}\,h(x)\mbox{,}\quad r\in\mathbb{N}\backslash\{0\},
				\end{equation*}
				\begin{equation*}\label{fbg2}
				{\bf E}\bigg(\bigg[\begin{array}{c} X  \\ r\end{array} \bigg]^Q_{p,q}\bigg)  = \sum_{x=r}^{\infty}\bigg[\begin{array}{c} x  \\ r\end{array} \bigg]^Q_{p,q}\,h(x)\mbox{,}\quad r\in\mathbb{N}\backslash\{0\}.
				\end{equation*}
				here designated  by  $r^{th}-$order generalized  $q-$ Quesne factorial and  $r^{th}-$order generalized $q-$ Quesne binomial moment, respectively, of the random variable $X.$
				In the particular case of $r=1,$ we deduce the  generalized $q-$ Quesne mean value, also called the generalized  $q-$ Quesne expectation value, of $X$  by 
				\begin{equation*}
				\mu^Q_{p,q}:={\bf E}\Big([X]^Q_{p,q}\Big)=\sum_{x=1}^{\infty}[x]^Q_{p,q}\,h(x).
				\end{equation*}
				The  associated variance of $X$ is then obtained as
				\begin{equation*}
				\big(\sigma^2_{p,q}\big)^Q:={\bf V}\Big([X]^Q_{p,q}\Big)= {\bf E}\big([X]^Q_{p,q}\big)^2 - \bigg[{\bf E}\Big([X]^Q_{p,q}\Big)\bigg]^2
				\end{equation*}
				or, equivalently,
				\begin{equation*}
				{\bf V}\Big([X]^Q_{p,q}\Big)=p^{-1}\,{\bf E}\big(\big({[X]^Q_{2,p,q}}\big)\big) +p^{X-2}q\,{\bf E}\big([X]^Q_{p,q}\big)- \big({q\over p}\big)^2\big[{\bf E}\big([X]^Q_{p,q}\big)\big]^2.
				\end{equation*}
				\item[(2)]
				Its deformed binomial moment is given 
				by
				\begin{equation*}\label{fbg6}
				{\bf E}\bigg(\bigg[\begin{array}{c} X  \\ j\end{array} \bigg]^Q_{p,q}\bigg)= \sum_{m=j}^{\infty} (-1)^{m-j}{(p-q^{-1})^{m-j}\over p^{-{m\choose 2}+\tau(m-j)}}s^Q_{p,q}(m,j){\bf E}\bigg(\bigg[\begin{array}{c} X  \\ m \end{array} \bigg]^Q_{p,q}\bigg),
				\end{equation*}
				while the factorial moment is expressed in terms of the  \textbf{generalized $q-$ Quesne} 
				factorial moment by
				\begin{equation*}\label{fbg7}
				{\bf E}[(X)_j]= j!\sum_{m=j}^{\infty} (-1)^{m-j}{(p-q^{-1})^{m-j}\over p^{-{m\choose 2}+\tau(m-j)}}s^Q_{p,q}(m,j){{\bf E}\big({[ X ]^Q_{m,p,q}}\big)\over [m]^Q_{p,q}!},
				\end{equation*}
				where $j\in\mathbb{N}\backslash\{0\},$ $\tau\in\mathbb{N},$ and $s^Q_{p,q}$ is the generalized $q-$ Quesne Stirling number of the first kind.
				\item[(3)] The probability distribution $h(x), x\in\mathbb{N},$ of a discrete random variable $X$ is given by
				\begin{equation*}\label{fbg9}
				h(x) = \sum_{m=x}^{\infty}(-1)^{m-x}\,p^{{x\choose 2}}\,q^{-{m-x\choose 2}}\bigg[\begin{array}{c} m \\ x \end{array} \bigg]^Q_{p,q}{\bf E}\bigg(\bigg[\begin{array}{c} X \\ m \end{array} \bigg]^Q_{p,q}\bigg),\quad x\in\mathbb{N}.
				\end{equation*}
			\end{enumerate}
		\end{remark}
		\section{Concluding remarks}
		In this paper, we have developed and illustrated the fundamentals of   $\mathcal{R}(p,q)-$ deformed combinatorics, with a special focus on factorials, binomial coefficients, Vandermonde's and Cauchy's formulae, binomial formula,  Stirling numbers, and Bell numbers induced by the $\mathcal{R}(p,q)-$ deformed quantum algebra. These results have also  been derived and discussed in the particular case of  the so-called generalized $q-$ Quesne deformed quantum algebra.
		Relevant properties have been deduced and analyzed in this framework. 
		\section*{Acknowledgements}
		This work is supported by TWAS Research Grant  RGA No.17 - 542 RG/ MATHS/AF/AC \_G - FR3240300147. 
		The ICMPA-UNESCO Chair is in partnership with the Association pour la 
		Promotion Scientifique de l'Afrique (APSA), France, and   Daniel Iagolnitzer Foundation (DIF), 
		France, supporting the development of mathematical physics in Africa. MNH acknowledges his colleagues Nicholas M. J. Hall,  Isabelle Dadou,   Yves Morel,  Catherine Jeandel, the staff of the UMR 5566 Laboratoire d'Etudes en G\'eophysique et Oc\'eanographique Spatiales (LEGOS) of the Facult\'e des Sciences et Ing\'enierie,  Michael Toplis, Director of the Observatoire Midi-Pyr\'en\'ees, and  Nguyen Tien Zung of the Institut de Math\'ematiques de Toulouse for their hospitality during his stay, as visiting professor,  at the Universit\'e Toulouse III Paul Sabatier, where this work has been completed. MNH is also grateful to Mrs Sophie Raynaud, Head of International Relations Office, and Prof. Fabrice Dumas, Vice-president in charge of International Relations,  Universit\'e Toulouse III Paul Sabatier,  for all their solicitude.


\begin{thebibliography}{15}
			\bibitem{A}  W. A. Al-Salam, $q-$ analogues of Cauchy's formulas, 
			{\it Proceedings of the American Mathematical Society}
			Vol. 17, No. 3 , pp. 616-621 (1966). 
			\bibitem{BS}Z. R. K. Balogh and M.J. Schlosser, $q-$ Stirling numbers of the second kind and $q-$ Bell numbers for graphs, {\it Elect. Notes. Disc. Maths.} {\bf 54}, 361-366 (2016).
			\bibitem{B}A. Broder, The $r-$ Stirling numbers, {\it Discrete Math.} {\bf 49}, 241-259 (1984).
			\bibitem{B1}I. M. Burban, Two-parameter deformation of oscillator algebra, {\it Phys. Lett.
				B.} {\bf 319}, 485-489 (1993).
			\bibitem{CK}W.S. Chung,  and H.J. Kang, $q-$ permutations and $q-$ combinations,  {\it Int. J. Theor. Phys.}{\bf 33}, 851-856  (1994) . 
			\bibitem{CJ} R. Chakrabarti and R. Jagannathan, {A $(p, q)-$oscillator realisation of two-parameter quantum algebras},  {\it  J. Phys.  A: Math. Gen.} {\bf 24}, L711-L718,  IMSC-91-15 (1991).
			\bibitem{CA1}Ch. A. Charalambides, {\it Discrete $q-$ distributions.}  John Wiley and Sons, Inc., Hoboken, New Jersey, 2016.
			\bibitem{CA2}Ch. A. Charalambides,  {\it Enumerative combinatorics.} Chapman and  Hall/CRC, Boca
			Raton, FL, 2002.
			\bibitem{CA3}Ch. A. Charalambides, Moments of a class of discrete  $q-$ distributions, 
			{\it J. Statist. Plann. Inference.} {\bf 135}, 67-85 (2004). 
			\bibitem{CA4}Ch. A. Charalambides,  {\it Combinatorial Methods in Discrete Distributions.} JohnWiley
			and Sons, Inc., Hoboken, NJ,2005.
			\bibitem{RC}R. B. Corcino, On $(p,q)-$ binomial coefficients, {\it E. J. Comb. Num. Theo.} {\bf 8},  (2008).
			\bibitem{RB} R. B. Corcino, C. Barrientos,  Some theorems of the $q-$ analogue of the generalized Stirling
			numbers, {\it Bull. Malays. Math. Sci. Soc}, {\bf 34} (3), 487-501 (2011).
			\bibitem{GW} H. W. Gould,  The $q-$Stirling numbers of the first and second kinds, {\it Duke Math. J.} {\bf 28},
			281-289 (1961).
			\bibitem{HB} M.N. Hounkonnou  and J. D. Kyemba Bukweli, $\mathcal{R}(p,q)$-calculus: differentiation and integration, {\it SUT. J. Math.} {\bf 49}, 145-167 (2013). 
			\bibitem{HB1}M. N. Hounkonnou  and J. D. Kyemba Bukweli, $(R,p,q)$-deformed quantum algebras: Coherent states and special functions, {\it J. Math. Phys.} \textbf{51} , 063518-063518.20 (2010).
			
			\bibitem{Hounkonnou&Ngompe07a} M. N. Hounkonnou and E. B. Ngompe Nkouankam, {New $(p, q, \mu, \nu, f)$-deformed states,}  {\it  J. Phys. A: Math. Theor.}  {\bf 40}, 12113-12130 (2007).
			\bibitem{JS}R. Jagannathan and K. Srinivasa Rao, Two-parameter quantum algebras,
			twin-basic numbers, and associated generalized hypergeometric series, {\it Proceedings of the International
				Conference on Number Theory and Mathematical Physics}, 20-21 December 2005.
			\bibitem{J1} F. H. Jackson, $q-$Difference equations, {\it Amer. J. Math.} {\bf 32}, 305-314 (1910).
			\bibitem{J1a}
			F. H. Jackson,  On $q-$definite integrals,  {\it Quart. J. Pure Appl. Math.} {\bf 41}, 193-203 (1910).
			\bibitem{KMM}G. Kalnins, W. Miller and S. Mukhejee, Models of $q-$algebra representations:
			matrix elements of the  $q-$ oscillator algebra, {\it J. Math. Phys}. {\bf 34}, 5333-
			5356 (1993).
			\bibitem{K1}T. Kim, A note on degenerate Stirling polynomials of second kind, Proc. Jangjeon Math.
			Soc. {\bf 20}, No. 3. pp. 319 - 331 (2017).
			\bibitem{K2}T. Kim, $\lambda-$ Analogue of Stirling numbers of the first kind, {\it Adv. Stud. Contemp. Math.} {\bf 27} (3), 423-429 (2017).
			\bibitem{K3}T. Kim, Y. Yao, D. S. Kim and G. W. Jang, Degenerate r-Stirling Numbers and $r-$ Bell Polynomials, {\it  Russ. J. Math. Phys.} {\bf 25},  44-58, (2018).
			\bibitem{L}C. L. Liu,    
			{\it Introduction to Combinatorial Mathematics}, Computer Science Series, New York, McGraw-Hill, 1968. 
			\bibitem{M}J. Mycielski, G.  Rozenberg and A. Salomaa,  
			{\it Structure in Logic and Computer Science}
			, Berlin, Springer-Verlag,  1997.  
			\bibitem{Q}C. Quesne, K. A. Penson and V. M. Tkachuk, Maths-type $q-$deformed coherent
			states for $q > 1.$ {\it Phys. Lett. A.} {\bf 313}, 29-36 (2003) .
			\bibitem{Q1}C. Quesne, New $q-$deformed coherent states with an explicitly known resolution
			of unity, {\it J. Phys. A: Math. Gen}. {\bf 35},   9213-9226 (2002).
			\bibitem{R}F. S.  Roberts,   
			{\it Applied Combinatorics}
			, Englewood Cliffs, New Jersey, Prentice-Hall, 1984.  
			\bibitem{WW} M. Wachs and D. White, $p,q-$ Stirling numbers and set partition statistics, {\it J. Comb. Theo, series A.
			} {\bf 56}, 27-46 (1991).
			
		\end{thebibliography}
	\end{document}